\documentclass[sn-mathphys-num]{sn-jnl}

\usepackage{graphicx}%
\usepackage{multirow}%
\usepackage{amsmath,amssymb,amsfonts}%
\usepackage{amsthm}%
\usepackage{mathrsfs}%
\usepackage[title]{appendix}%
\usepackage{xcolor}%
\usepackage{textcomp}%
\usepackage{manyfoot}%
\usepackage{booktabs}%
\usepackage{algorithm}%

\usepackage{algpseudocode}%
\usepackage{listings}%
\usepackage{microtype}

\usepackage{bm}
\usepackage{parskip}
\usepackage{cleveref}
\crefname{algorithm}{Algorithm}{Algorithms}

\crefname{thm}{Theorem}{Theorems}
\crefname{theorem}{Theorem}{Theorems}
\crefname{lemma}{Lemma}{Lemmas}
\crefname{equation}{Eq.}{Eqs.} 
\crefname{figure}{Figure}{Figures} 
\crefname{section}{Section}{Sections}
\crefname{remark}{Remark}{Remarks}
\Crefname{algorithm}{Algorithm}{Algorithms}

\usepackage{mathtools}
\usepackage{float}
\usepackage[makeroom]{cancel}
\theoremstyle{thmstyleone}%
\newtheorem{theorem}{Theorem}
\newtheorem{thm}[theorem]{Theorem}
\newtheorem{lemma}[theorem]{Lemma}
\theoremstyle{thmstyletwo}%

\newtheorem{remark}{Remark}%

\theoremstyle{thmstylethree}%
\newtheorem{definition}{Definition}%
\numberwithin{remark}{section}
\numberwithin{definition}{section}
\numberwithin{equation}{section}
\numberwithin{theorem}{section}
\usepackage{csquotes}
\usepackage{cleveref}
\usepackage{subcaption}
\usepackage{mathrsfs}
\usepackage{xtab}

\newcommand{\norm}[1]{\left\lVert#1\right\rVert}

 \newcommand{\N}{\ensuremath{\mathbb{N}}}
\newcommand{\R}{\ensuremath{\mathbb{R}}}

\newcommand{\tit}[1]{\textit{#1}}
\begin{document}

\title[]{Sparse Randomised Approximation of Normal Cycles}

\author[1]{\fnm{Allen} \sur{Paul}}\email{allen.paul.alb@gmail.com}

\author[2]{\fnm{Neill} \sur{Campbell}}\email{nc537@bath.ac.uk}

\author*[1]{\fnm{Tony} \sur{Shardlow}}\email{t.shardlow@bath.ac.uk}

\affil*[1]{\orgdiv{Department of Mathematical Sciences}, \orgname{University of Bath}, \orgaddress{ \city{Bath}, \postcode{BA2 7AY},  \country{UK}}}
\affil[2]{\orgdiv{Department of Computer Science}, \orgname{University of Bath}, \orgaddress{ \city{Bath}, \postcode{BA2 7AY},  \country{UK}}}

\date{}

\abstract{We extend our work~\cite{SparseNystromCurrVar} for compression of currents and varifolds to a compression algorithm for the embedded normal cycles representation of shape, restricted to the constant normal kernel case, using the Nystrom approximation in Reproducing Kernel Hilbert Spaces (RKHS) and Ridge Leverage Score (RLS) sampling. Our method comes with theoretical guarantees on the  compression error decay, and the approximations are shown to be effective for downstream tasks such as nonlinear shape registration in the Large Deformation Diffeomorphic Metric Mapping (LDDMM) framework, even for very high compression ratios. The performance of our algorithm is demonstrated on large-scale shape data from modern geometry processing datasets and is shown to  accelerate  downstream registration tasks significantly.}
\keywords{normal cycles, compression, Nystrom approximation, randomised approximation, shape registration, LDDMM}
\maketitle

\section{Introduction}\label{intro}
In many real-world geometric learning applications --- for example, in computational anatomy --- one wishes to build statistical models of shape variation, where the shape data are available as discrete curves and surfaces representing anatomical structures~\cite{Younes}. In order to fit such models, one requires a fidelity metric on shapes in order to compare the data to shapes generated by the model and tune the model parameters. In the majority of such applications, shapes do not have a consistent parametrisation, meaning that simple fidelity metrics relying on known correspondences  cannot be used. A large literature has been developed on point-cloud metrics~\cite{Hausdorffref,Chamferref} such as Chamfer and Hausdorff distance, as well as methods for extracting dense correspondences~\cite{FINDINGCORRESPONDENCE} between shapes. However, in general, point-cloud distances are \tit{not} geometric and do not account for the underlying structure of the shapes. Furthermore, dense correspondences may not be well defined or unique, which introduces a bias in subsequent shape learning applications.

In the case where shape data are available as submanifolds of $\R^{d}$, one can deal with the lack of parametric correspondences in a more principled manner, using techniques from geometric measure theory \cite{FidelityMetrics}. In particular,  the so-called currents \cite{GlaunesCurrents} and varifolds~\cite{Charon} representations  view shapes as objects that integrate a family of continuous differential forms on the underlying domain, effectively embedding the shapes into the dual space of a space of differential forms. By restricting to differential forms lying in a Reproducing Kernel Hilbert space (RKHS) \cite{RKHS}, one may use the dual metric on the embeddings in order to compare shapes in terms of their action on forms, independent of their parametrisations. The computation of the dual metric between submanifolds can be written down in closed form in terms of the positive semidefinite kernel function $k\colon\R^{d}\times \R^{d} \to \R$ associated to the RKHS. One may compute the RKHS embeddings and metrics between discrete shape data in a way that is consistent as the discrete shapes tend to the underlying continuous submanifold. This framework and the resulting shape metrics  have been used extensively  in the Large Deformation Diffeomorphic Metric Mapping (LDDMM) framework~\cite{FidelityMetrics}  for matching shapes with diffeomorphisms.

While currents and varifolds metrics allow for \tit{geometric} comparison of shapes in a correspondence-less manner, these representations only take into account first-order geometric information, such as tangent and normal vector directions. The structures  typically present in real-world applications have regions of high curvature, branching points or distinguished boundaries; in such cases, currents and varifolds metrics are limited as they do  not capture `higher-order' properties of the underlying submanifold. Indeed,  undesirable features may be generated during diffeomorphic registration when these metrics are used as discrepancy terms for matching complex shapes with the aforementioned features~\cite{FidelityMetrics}. 

A theoretically grounded solution is the normal cycles representation, which originates in geometric measure theory~\cite{Federer,Federer1959} and has been applied to computational anatomy and LDDMM~\cite{NCThesis,NCpaper}. Whereas the current associated to a surface integrates test differential forms over the surface itself, the normal cycle integrates them over the surface's \tit{unit normal bundle}, yielding an enriched representation that is sensitive not only to the surface but also to the rate of change of its normal vectors. As a consequence, the normal cycle locally encodes the mean and Gaussian curvatures of the surface~\cite{NCThesis}.

In a similar manner to currents and varifolds, one may compute the distance between two shapes in terms of the dual metric between their normal cycles, embedded into the dual of a suitable RKHS of forms. The advantages of comparing shapes with normal cycles instead of currents or varifolds are evidenced numerically  on shapes with high curvature regions, boundaries and branching points in~\cite{NCpaper}. However,  the practical computation of the normal cycles metric is costly in terms of both  memory and computational complexity. The normal cycles metric complexity scales as $\mathcal{O}(MN)$ when comparing triangulated surfaces with $M,N$ edges, respectively. While tools such as the KeOps library \cite{Keops} exist for large-scale metric (and gradients) computation, it is observed empirically in \cite{NCpaper} that even with these tools, the normal cycles representation is at least $3$--$5$ times more expensive to compute and take gradients of than currents and varifolds and can become very slow for large-scale problems with resolution past $10^{5}$. 


\subsection{Contribution}
In this work, we extend the randomised projection algorithm of \cite{SparseNystromCurrVar} from currents and varifolds to normal cycles. Building on \cite{NCpaper}, we derive an explicit Dirac delta basis form for the discrete embedded normal cycle in the constant normal kernel case (as in \cite{NCpaper}), which lets us apply the main theorem of \cite{SparseNystromCurrVar} and compress a target normal cycle of size $M$ to a sparse approximation of size $m\ll M$. Post-compression, the metric between two compressed normal cycles can be computed in $\mathcal{O}(m^{2})$ rather than $\mathcal{O}(M^{2})$, with the corresponding gain in gradient cost. To our knowledge this is the first compression algorithm for normal cycles, allowing the metric to scale to problems with resolution $M\sim 10^{5}$. The theoretical guarantees of \cite{SparseNystromCurrVar} carry over directly, yielding exponential error-decay bounds (for Gaussian RKHS kernels) of the form
\begin{align*}
    \norm{\mu - \hat{\mu}}_{W^{*}}^{2} \leq \mathcal{O}(m\exp(-\alpha m^{\frac{1}{d}})),
\end{align*}
where $\hat{\mu}$ is the approximation of size $m\ll M$ to the original normal cycle $\mu \in W^{*}$ and $\alpha>0$. The practical benefits are demonstrated in~\cref{numexp}, where we obtain order-of-magnitude acceleration of LDDMM registration on heavily over-sampled real-world surface data; the downstream utility of \cref{algorithm1} is further demonstrated in the companion work~\cite{CurvatureEnthusiasm}, where it is the core matching metric in a learning-based pipeline.

\subsection{Existing work}
Compression of currents and varifolds has been studied previously in~\cite{DURRLEMAN,VarCompression}, but the methods of those works are substantially slower than the ridge-leverage-score scheme of~\cite{SparseNystromCurrVar} on which we build, and do not readily extend to the normal cycles setting. Other work on normal cycles~\cite{NCThesis,NCpaper,CohenSteiner2003} has focused on formulation and on applications to surface registration and curvature estimation, not on scaling the metric to large $M$.

\subsection{Outline}
In \cref{bg}, we review the currents and normal cycles representations of shape within the framework of exterior algebra, differential forms and rectifiable sets. In \cref{mainsec}, we describe how, in the constant normal kernel setting, one may derive an explicit Dirac delta decomposition of the embedded normal cycle of a triangulated surface, and how one may apply the randomised compression algorithm of \cite{SparseNystromCurrVar}. Finally, in \cref{numexp}, we demonstrate the strengths and weaknesses of the proposed compression algorithm on large-scale shape data from modern geometry processing datasets.

\section{Background}\label{bg}
The theory of currents was developed by Federer~\cite{Federer1959}
within the framework of geometric measure theory, originally for the
study of the Plateau problem; the notions of curvature measure,
sets of positive reach, and the additivity property of normal
cycles originate in the same work, see also the
monograph~\cite{Federer}. The normal cycle, which generalises 
curvature measures to non-smooth sets, was subsequently developed 
by Wintgen~\cite{wintgen1982normal}, 
Z\"ahle~\cite{zahle1986integral}, and 
Fu~\cite{fu1989monge}, among others. Curvature estimation on 
polyhedral surfaces via normal cycles was studied by Cohen-Steiner 
and Morvan~\cite{CohenSteiner2003}; see 
also~\cite{Morvan2008} for a comprehensive treatment of 
generalised curvatures and normal cycles. The application of normal 
cycles to computational anatomy and LDDMM was introduced 
in~\cite{NCThesis,NCpaper}.

We now define the normal cycle rigorously, using exterior algebra
and differential forms; see \cref{Appendix1} for background. We
follow the notation of~\cite{NCThesis} and begin with the currents
representation, from which the normal cycle is obtained as a
higher-order extension.
\subsection{Currents}
The core idea of currents is to represent submanifolds as objects that integrate continuous differential forms. Formally, one begins by defining the space of continuous compactly supported differential $m$-forms (for $m\leq d$)
\begin{align*}
    \Omega_{0}^{m}(\R^{d}) = \mathcal{C}_{0}(\R^{d},(\Lambda^{m}\R^{d})^{*}),
\end{align*}
equipped with the uniform norm, where $(\Lambda^{m}\R^{d})^{*}\simeq \Lambda^{m}\R^{d}$ is the dual to the space of $m$-vectors formed through wedge products. See \cref{Appendix1} for background on exterior product spaces. This leads to the following definition.
\begin{definition}
The space of $m$-currents is defined as $\Omega_{0}^{m}(\R^{d})^{*}$, which is the topological dual of $ \Omega_{0}^{m}(\R^{d})$. This means that $m$-currents are bounded linear functionals $T\colon \Omega_{0}^{m}(\R^{d}) \to \R $ such that the following holds for some $C(T)>0$:
\begin{align*}
    \lvert T(\omega)\rvert \leq C(T)\norm{\omega}_{\infty},\qquad \forall \omega \in  \Omega_{0}^{m}(\R^{d}).
\end{align*}
Given an (oriented) $m$-rectifiable set $S\subset \R^{d}$, the unique current associated with it is defined as $[S]\in  \Omega_{0}^{m}(\R^{d})^{*} $ with action
\begin{align}\label{currdef}
    [S](\omega) \coloneq \int_{S}{ (\omega(x)|\tau^{S}(x) ) \,d\mathcal{H}^{m}(x)},\quad \omega \in \Omega_{0}^{m}(\R^{d}),
\end{align}
where for $\mathcal{H}^{m}$-almost every $x$ in $S$, one defines $\tau^{S}(x) \coloneq e_{1}(x)\wedge \dots \wedge e_{m}(x)$, a wedge product of (positively oriented) orthonormal basis vectors of the tangent space at $x$. This is well defined, as the tangent space exists $\mathcal{H}^{m}$-almost everywhere for an $m$-rectifiable set, and the choice of basis vectors does not matter due to invariance of the wedge product under change of orthonormal basis. 
\end{definition}

From \cref{currdef}, we observe that the current associated to $S$  integrates a differential $m$-form over  $S$ with respect to the Hausdorff measure. When $S$ is a smooth oriented $m$-dimensional submanifold, this reduces to the standard definition of integration of $m$-forms against $S$ \cite{NCThesis}. As shown in \cref{simplelemmma},
\begin{align*}
    \lvert [S](\omega)\rvert \leq \norm{\omega}_{\infty}\mathcal{H}^{m}(S),\quad \omega \in  \Omega_{0}^{m}(\R^{d})
\end{align*}
so that $[S]$ is indeed a well-defined element of $ \Omega_{0}^{m}(\R^{d})^{*}$. One can also define a natural push-forward action, $\varphi_{\#}[S] \in \Omega_{0}^{m}(\R^{d})^{*}$, on currents under the action of a diffeomorphism $\varphi \colon \R^{d} \to \R^{d}$ on the ambient space. While we do not pursue this definition here, one may prove \cite{NCThesis} that this action is geometric in the sense that
$   \varphi_{\#}[S] = [\varphi(S)]$.
This identity makes the diffeomorphic pushforward of a current a straightforward object to handle. This is particularly useful for applications of currents in the LDDMM framework, where shapes are modelled as deformations of template shapes under the action of a diffeomorphism.

\subsubsection{Dual space metrics for currents}
It was originally suggested in \cite{GlaunesCurrents} that in practical applications, one may compare shapes and submanifolds in a non-parametric, correspondence-less manner, by comparing the associated currents in the dual space metric. The associated distance between the currents associated to shapes may then be measured as
\begin{align}\label{currdualmet}
   d\bigl([S],[T]\bigr) \coloneq \norm{\vphantom{\big|}[S]-[T]}_{\Omega_{0}^{m}(\R^{d})^{*} } = \underset{\substack{\omega \in \Omega_{0}^{m}(\R^{d})\\\norm{\omega}_\infty \le 1}}{\sup}{\big\lvert [S](\omega) - [T](\omega)\big\rvert }.
\end{align}
In \cref{currdualmet}, shapes are compared in terms of how similarly or differently they integrate the same differential forms, and one computes the maximum such discrepancy over all $\omega \in \Omega_{0}^{m}(\R^{d})$. Unfortunately, the uniform norm topology over $\Omega_{0}^{m}(\R^{d})$ is too strong and does not give a useful measure of dissimilarity to compare and register shapes as proven in \cite{Charon}. Indeed, the distance \cref{currdualmet} simply collapses to the sum of the Hausdorff measure of each shape when $S\neq T$ -- thus  failing to give a useful distance between the shapes.  


Instead, as proposed by \cite{GlaunesCurrents}, one embeds currents of interest $[S], [T]$ into the \textit{dual} of a smaller space $W \subset \Omega_{0}^{m}(\R^{d})$ of differential forms, equipped with an RKHS structure. The natural dual metric on $W^{*}$ induces a weaker topology on currents, and has closed-form expressions which one may use for comparison and optimisation purposes. Following~\cite{GlaunesCurrents}, we let $W$ be an RKHS of differential forms, generated by a kernel $K\colon \R^{d} \times \R^{d} \to \mathcal{L}(\Lambda^{m}\R^{d})$ of the form
\begin{align*}
    K(x,y) = K_{p}(x,y)\mathrm{Id}_{\Lambda^{m}\R^{d}},
\end{align*}
where $K_{p}$ is a sufficiently smooth positive-definite scalar-valued spatial kernel on $\R^{d}$. Such a $W$ takes the form
\begin{align*}
    W = \bigg\{f= \sum_{i=1}^{\infty}{K_{p}(\cdot,x_{i})\tau_{i} } \colon\tau_{i} \in \Lambda^{m}\R^{d} ,\quad x_{i} \in \R^{d},\quad \norm{f}_{W}<\infty\bigg\},
\end{align*}
where $\norm{f}_{W}^{2}=\sum_{i,j=1}^{\infty}{K_{p}(x_{i},x_{j})\langle \tau_{i},\tau_{j}\rangle_{\Lambda^{m}\R^{d}} }$.
Provided that $W$ is dense in $\Omega_{0}^{m}(\R^{d})$ and $K_p$ is bounded, this guarantees the dual embedding $\Omega_{0}^{m}(\R^{d})^{*} \subset W^{*}$ is an injection \cite{NCThesis}, and that the $W^{*}$ metric defines a true metric on the space of currents. A popular choice of kernel function that ensures this is the Gaussian radial basis function (RBF) kernel. The action of $[S],[T]$ as elements of $W^{*}$ is the same integral action as in \cref{currdef}, and they are well-defined elements of $W^{*}$ due to the bound
\begin{align*}
    \lvert [S](\omega)\rvert \leq \norm{\omega}_{\infty}\mathcal{H}^{m}(S) \leq c\norm{\omega}_{W}\mathcal{H}^{m}(S),\qquad \omega \in W,
\end{align*}
where we have used the RKHS identity $\norm{\omega}_{\infty} \leq c \norm{\omega}_{W}$ for some $c>0$. Using the reproducing property of the RKHS kernel, one may compute \cite{GlaunesCurrents} the dual metric in $W^{*}$ between $[S],[T]$ explicitly in terms of integration as
\begin{align}\label{dualmet_w}
    d_{W^{*}}([S],[T])^{2} = \norm{\vphantom{\big|}[S]-[T]}_{W^{*}}^{2} = \bigl\langle [S],[S] \bigr\rangle_{W^{*}} -2\bigl\langle [S],[T] \bigr\rangle_{W^{*}} + \bigl\langle [T],[T] \bigr\rangle_{W^{*}},
\end{align}
where the product terms are defined as
\begin{align}\label{intergral_w}
    \bigl\langle [S],[T] \bigr\rangle_{W^{*}} =  \int_{S}\int_{T}{K_{p}(x,y)\langle \tau^{S}(x),\tau^{T}(y) \rangle_{\Lambda^{m}\R^{d} }\,d\mathcal{H}^{m}(x) \,d\mathcal{H}^{m}(y)}.
\end{align}
In the discrete case, when we have triangular approximations $\hat{S},\hat{T}$ of  surfaces $S,T$ for example, one forms approximate currents
\begin{align}\label{currentdeltaapprox}
    [\hat{S}] \approx \widehat{[S]} \coloneq \sum_{i=1}^{n}{\delta_{c_{i,\hat{S}} }\tilde{\tau}_{i,\hat{S}} } \in W^{*}, \quad
    [\hat{T}] \approx \widehat{[T]} \coloneq \sum_{i=1}^{m}{\delta_{c_{i,\hat{T}}}\tilde{\tau}_{i,\hat{T}}  } \in W^{*},
\end{align}  
where $c_{i,\hat{S}},c_{i,\hat{T}}$ denote the centre of the $i$th triangle of $\hat{S},\hat{T}$, respectively. The discrete weights are defined as $\tilde{\tau}_{i,\hat{S}}={\tau}^{\hat{S}}(c_{i,\hat{S}})a_{i,S}$ where $a_{i,\hat{S}}$ is the $i$th triangle area. Applying the reproducing property of the RKHS kernel yields an easily computed dual metric  with `discrete' product term
\begin{align}\label{discinnerprodcurr}
    \big\langle \widehat{[S]},\widehat{[T]}\big\rangle_{W^{*}} = \sum_{i=1}^{n}\sum_{j=1}^{m}{K_{p}(c_{i,\hat{S}},c_{j,\hat{T}} )\langle \tilde{\tau}_{i,\hat{S}},\tilde{\tau}_{j,\hat{T}} \rangle_{\Lambda^{m}\R^{d}} }.
\end{align}
This may be interpreted as the discretisation of \cref{intergral_w}. This yields the discrete currents metric
\begin{align*}
    d_{W^{*}}\bigl(\widehat{[S]},\widehat{[T]}\bigr)^{2} 
    = \bigl\|{\widehat{[S]}-\widehat{[T]}}\bigr\|^2_{W^{*}} = \bigl\langle \widehat{[S]},\widehat{[S]} \bigr\rangle_{W^{*}} -2\bigl\langle \widehat{[S]},\widehat{[T]} \bigr\rangle_{W^{*}} + \bigl\langle \widehat{[T]},\widehat{[T]} \bigr\rangle_{W^{*}}.
\end{align*}
Such delta approximations as in \cref{currentdeltaapprox} are known to be consistent \cite{FidelityMetrics} in the sense that
$    \bigl\|{\widehat{[S]} -  [S]  }\bigr\|_{W^{*}} \lesssim \tau(h),
$ 
where the bounding function satisfies $\tau(h) \to 0$ as $h$, the triangulation mesh size, tends to zero.

In order to practically compute the weights in \cref{discinnerprodcurr}, one observes \cite{FidelityMetrics} the identification $\Lambda^{2}\R^{3} \simeq \Lambda^{1}\R^{3} \simeq \R^{3}$. For surface data, we have $(d,m)=(3,2)$ so that we may identify $W$ with a space of vector fields on $\R^{3}$. Under this identification, for triangulated surfaces, $\tau^{S}(x)$ corresponds to the unit normal vector at $x \in S$, which is simple to compute. For curves in $\R^{3}$ or $\R^{2}$, we have $(d,m)=(3,1)$ or $(d,m)=(2,1)$, respectively, so $W$ is again a space of vector fields. In either case, similar Dirac delta decompositions hold, with the weights reducing to the discrete tangent vectors of the discrete curve. 

\subsection{Normal cycles}
The normal cycle of a smooth oriented submanifold $S\subset \R^{d}$ is the current associated to the unit normal bundle of $S$; integrating against the normal bundle captures both position and the rate of change of normals along $S$. We give the precise definition following the notation of~\cite{NCpaper}.
First, we introduce the sets of positive reach following~\cite{Federer1959}.
\begin{definition}
    If $S\subset \R^{d}$, the reach $R \in \R$ of $S$ is defined as the supremum of $r>0$ such that there exists a unique projection onto $S$ from $\partial S_{r}\coloneq\{x \in \R^{d}: d(x,S)=r\}$, where $d(x,S) \coloneq \inf_{s\in S}\|x-s\|$. A set with reach $R>0$ is a set with positive reach.
\end{definition}
This definition encompasses all convex sets and smooth submanifolds. For any such $S$ and any $r\in(0,R)$, the unit normal bundle is
\begin{align}\label{unitnormalbunfle}
    \mathcal{N}_{S} \coloneq \{(x,n) \in S \times \mathbb{S}^{d-1}: x+rn \in \partial S_{r} \},
\end{align}
which is independent of $r$ and of the orientation of $S$, and is a $(d-1)$-dimensional orientable rectifiable subset of $\R^{d}\times \mathbb{S}^{d-1}$. For a smooth surface or space curve $S\subset \R^{3}$, for example, $\mathcal{N}_{S}$ is a $2$-dimensional manifold in $\R^{3}\times \mathbb{S}^{2}$, visualised as a smooth `fattening' of $S$. Rectifiability and canonical orientation make the associated current well defined; this current is the \tit{normal cycle} of $S$, defined formally below.
Let
 \begin{align*}
 \Omega_{0}^{d-1}(\R^{d}\times \mathbb{S}^{d-1}) = \mathcal{C}_{0}(\R^{d}\times \mathbb{S}^{d-1},(\Lambda^{d-1}(\R^{d}\times \R^{d} )) ^{*})
 \end{align*}
 for $(d-1)$ differential forms on $\R^{d}\times \mathbb{S}^{d-1}$.
\begin{definition}
    Given a set of positive reach $S\subset \R^{d}$, the normal cycle associated to $S$ is defined as $N(S)\coloneq[\mathcal{N}_{S}]$ which is an element of $\Omega_{0}^{d-1}(\R^{d}\times \mathbb{S}^{d-1})^{*}$. The action of $N(S)$ on $\omega \in \Omega_{0}^{d-1}(\R^{d}\times \mathbb{S}^{d-1})$ is given by
    \begin{align}\label{NCaction}
        N(S)(\omega) = [\mathcal{N_{S}}](\omega) = \int_{\mathcal{N}_{S}}{\Bigl(\omega(x,\nu) \mid \tau^{\mathcal{N_{S}}}(x,\nu)\Bigr)\,d\mathcal{H}^{d-1}(x,\nu)},
    \end{align}
    where $\tau^{\mathcal{N_{S}}}(x,\nu) = e_{1}(x,\nu) \wedge \dots \wedge e_{d-1}(x,\nu)$ is the wedge product of a positively oriented orthonormal basis of the tangent space at point $(x,\nu)$.
\end{definition}
Specific differential forms recover the mean and Gaussian curvatures of $S$ locally~\cite{NCThesis}.

In practical geometric learning, we work with discrete shapes in the form of the union of curve segments and triangulated surfaces. The definition of normal cycle given above may also be generalised \cite{NCpaper,NCThesis} to all finite unions of sets of positive reach. Indeed, if $S = S_{1}\cup S_{2}$ is a union of sets of positive reach, such that $S_{1}\cap S_{2}$ has positive reach, then
the additivity property holds \cite{Federer1959}:
\begin{align}\label{unionformula}
    N(S) \coloneq N({S_{1}}) + N({S_{2}}) - N({S_{1}\cap S_{2}}).
\end{align}
The left and right-hand sides agree if $S$ has positive reach. In order to compute the normal cycle $N(\mathcal{T})$ of a triangulated surface $\mathcal{T} = \cup_{i=1}^{n_{T}}{T_{i}}$ given as a union of closed triangles, one may apply \cref{unionformula} recursively to obtain a decomposition in terms of the normal cycles of individual triangles and segments. We shall use this decomposition in \cref{mainsec}.

\subsubsection{Dual space metric: normal cycles}
Similar to the case of currents, the topology on the dual space induced by the uniform norm is too strong to give a meaningful distance between normal cycles. As for currents, fixing an appropriate RKHS $W\subset \Omega_{0}^{d-1}(\R^{d}\times \mathbb{S}^{d-1})$ of test differential forms, one embeds  the normal cycle representation into the dual $W^{*}$. As we shall see, this allows one to tractably compute metrics  between normal cycles in terms of the kernel function.

One begins by defining a $\mathcal{L}(\Lambda^{d-1}(\R^{d}\times  \R^{d}))$-valued positive-definite kernel on $\R^{d} \times \mathbb{S}^{d-1}$ as
\begin{align}\label{productkernel}
    K((x,u),(y,v)) = K_{p}(x,y)K_{s}(u,v)\mathrm{Id}_{\Lambda^{d-1}(\R^{d}\times  \R^{d})},\quad \forall (x,u),(y,v)\in \R^{d}\times \mathbb{S}^{d-1}.
\end{align}
For sufficient regularity of $K_{p}$ and $K_{s}$, this induces an RKHS $W\subset \Omega_{0}^{d-1}(\R^{d}\times \mathbb{S}^{d-1})$ of forms $\omega \colon \R^{d}\times \mathbb{S}^{d-1} \to \Lambda^{d-1}(\R^{d}\times  \R^{d})$ (see~\cite{NCpaper, NCThesis}). The space $W$ takes the form
\begin{equation}\label{embeddingspace}
W = \bigg\{f = \sum_{i=1}^{\infty}{K_{p}(\cdot,x_{i})K_{s}(\cdot,s_{i})\tau_{i} }\colon \begin{subarray}{l}(x_{i},s_{i},\tau_{i}) \in \R^{d}\times \mathbb{S}^{d-1}\times \Lambda^{d-1}(\R^{d}\times  \R^{d}),\\ \norm{f}_{W}<\infty\end{subarray}\bigg\}
\end{equation}
where
\[
\norm{f}_{W}^{2}\coloneqq\sum_{i,j=1}^{\infty}{K_{p}(x_{i},x_{j})K_{s}(s_{i},s_{j})\langle \tau_{i},\tau_{j}\rangle_{\Lambda^{d-1}(\R^{d}\times  \R^{d}) } }.
\]
Typical choices for the spherical kernel are constant $K_{s}(u,v) = 1$, linear $K_{s}(u,v) = \langle u,v\rangle$, and Gaussian $K_{s}(u,v) = \exp(-\frac{\norm{u-v}^{2}}{2\sigma^{2}})$. One typically sets $K_{p}$ to be the Gaussian RBF kernel. 

With the same action as \cref{NCaction} on elements of $W$, a normal cycle $N(S)$ may be embedded as a well-defined element of $W^{*}$, through the inclusion $\Omega_{0}^{d-1}(\R^{d}\times \mathbb{S}^{d-1})^{*} \subset W^{*}$. By the reproducing property of the kernel, the dual metric between normal cycles $N(S),N(T)$ may be expressed  as
\begin{align*}
    \norm{N(S) - N(T) }_{W^{*}}^{2} = \langle N(S),N(S)\rangle_{{W^{*}}} - 2\langle N(S),N(T)\rangle_{_{W^{*}}} + \langle N(T),N(T) \rangle_{{W^{*}}},
\end{align*}
with product term
\begin{gather}\label{NCinnerprod}
\begin{split}
    \langle N(S),N(T)\rangle_{_{W^{*}}} \coloneq& \int_{\mathcal{N}_{S}}
    \int_{\mathcal{N}_{T}}
    K_{p}(x,y)
    K_{s}(u,v)\\
    &\qquad\qquad\quad
    \times
    \langle
    \tau_{S}(x,u),\tau_{T}(y,v)\rangle
    \,d\mathcal{H}^{d-1}(x,u) \,d\mathcal{H}^{d-1}(y,v),
\end{split}
\end{gather}
where $\tau_{S}(x,u),\tau_{T}(x,u)$ are defined as in \cref{NCaction}. Under suitable conditions on the spatial and spherical kernels,  the inclusion $\Omega_{0}^{d-1}(\R^{d}\times \mathbb{S}^{d-1})^{*} \subset W^{*}$ is injective, and the dual metric yields a well-defined distance on $\Omega_{0}^{d-1}(\R^{d}\times \mathbb{S}^{d-1})^{*}$. For the $K_{s}=1$ case, which is the primary setting for surface matching explored in \cite{NCpaper} and in this work, the resulting $W$ is not necessarily dense in the space of differential forms, and thus the $W^{*}$ metric is only a pseudo-metric on $\Omega_{0}^{d-1}(\R^{d}\times \mathbb{S}^{d-1})^{*}$. 

For practical applications with discrete shapes,  the inner-product in \cref{NCinnerprod}  must be approximated; this is significantly trickier than for currents, as the normal cycle also involves integration over the `spherical' components of the normal bundle. For triangulated surfaces in the case of $K_{s} = 1$, this is the content of the following theorem from \cite{NCpaper}.
\begin{thm} \label{thm:NCmetric}
    Let $\mathcal{T},\mathcal{T}'$ be two triangulated surfaces. The inner product between the associated discrete normal cycles can be computed in $W^{*}$ when $K_{s}=1$ as
    \begin{gather}\begin{split}\label{NC_disc_met}
        \langle &N(\mathcal{T}),N(\mathcal{T}') \rangle_{W^{*}} \\&=
        4\sum_{i=1}^{n_{e}}\sum_{j=1}^{m_{e}}{K_{p}(c_{i},d_{j})\langle f_{i},g_{j}\rangle}
        \times\bigg\langle \sum_{\{T\colon f_{i}\in \partial T\}}{n_{T,f_{i}}}, \sum_{\{T\colon g_{j}\in \partial T\}}{n_{T,g_{j}}}\bigg\rangle 
        \\&\quad+ \ \pi^2\sum_{x_{i}\in\partial \mathcal{T}}\sum_{y_j\in \partial \mathcal{T}'}{K_{p}(x_{i},y_{j})\langle A_{i},B_{j}\rangle},
    \end{split}\end{gather}
where $\partial \mathcal{T}$ denotes the boundary of the surface $\mathcal{T}$ and similarly for $\mathcal{T}'$. In the above, $n_{e},m_{e}$ denote the number of unique edges in each triangulation, and $f_{i},g_{j}$ the edges of $\mathcal T, \mathcal T'$ respectively. The kernel function in the first double sum is evaluated on $c_{i},d_{j}$, which denote centres of the edges $f_{i},g_{j}$, respectively. In the second double sum, the kernel function is evaluated on $x_{i},y_{j}$, which are the vertices of the boundary. Furthermore, the quantity $A_{i} \coloneq \sum_{k}{f_{k}^{i}/\lvert f_{k}^{i}\rvert}$ is defined as the sum of normalised edges $f^i_k$ attached to the boundary vertex $x_{i}$, oriented outwards from $x_{i}$. Finally, $n_{T,f_{i}}$ denotes the normal to triangle $T$ coinciding with edge $f_{i}$, such that $n_{T,f_{i}}\times f_{i}$ is oriented inwards for the triangle $T$.
 \end{thm}

\begin{remark}\label{rem:coefficient_discrepancy}
The constant prefactors $4$ and $\pi^{2}$ in \cref{NC_disc_met} differ from the single prefactor $\pi^{2}/4$ appearing for both terms in \cite[Prop.~3]{NCpaper}. Tracing the derivation back to the underlying integrals, this discrepancy originates in the evaluation of the half-circle integral $\int_{C}u\,d\mathcal{H}^{1}(u)$ over the semicircle $C=S^{\perp,+}_{f,-f\times n_{T}}$. A self-contained derivation is given in \cref{halfcircleint}. The formulae in \cref{thm:NCmetric,NCCURVEMET} use the corrected value, and the derivation of the Dirac delta decomposition in \cref{APPENDIXSECPROOF} is consistent with it.
\end{remark}

 Unlike currents and varifolds, the normal cycle metric has an explicit term for comparing boundaries of shapes, reflecting its sensitivity to boundary geometry.  The computational complexity of computing \cref{NC_disc_met} is much higher than that of currents. The first term is a double loop over edges, of which there are more than the number of triangles, and each term involves more inner product computations than for currents. Furthermore, the added boundary term also contributes an additional cost for shapes with boundary. Even with KeOps, the practical cost of computing the metric and its gradient is around $3$--$5$ times higher for normal cycles compared to currents and varifolds \cite{NCpaper}.  

For discrete curves in the case $K_s=1$, we have the following representation of the metric as a double sum across edges~\cite{NCThesis}.
\begin{thm}
    Let $C$ and $S$ be two discrete curves (unions of segments) with $N$ and $M$ vertices, respectively. Denoting the spatial kernel as $K_{p}$ and fixing constant normal kernel $K_{s}=1$, the inner product between the associated normal cycles $N(C), N(S)$ in $W^{*}$ is given by
    \begin{align}\label{NCCURVEMET}
        \langle N(C), N(S) \rangle_{W^{*}} = {\pi^{2}}{}{\sum_{i=1}^{N}\sum_{j=1}^{M}{K_{p}(x_{i},y_{j})}\langle A_{i},B_{j}\rangle }, 
    \end{align}
    where $A_{i} = \sum_{k}{f_{k}^{i}/\lvert f_{k}^{i}\rvert}$ is the sum of normalised edges attached to vertex $x_{i}$, and oriented outwards from $x_{i}$.
\end{thm}
 This has a similar complexity to the currents metric between discrete curves.

\begin{remark}
When $K_{s}=1$ as in the previous theorems, the product kernel  reduces to 
\begin{align*}
    K((x,u),(y,v) )=K_{p}(x,y) \mathrm{Id}_{\Lambda^{d-1}(\R^{d}\times  \R^{d})}
\end{align*}
independent of normal component, and elements of $W$ are constant on $\mathbb{S}^{d-1}$. Thus, one may isometrically identify $W$ with the space 
\begin{equation}\label{newembeddingspace}
 W = \bigg\{f = \sum_{i=1}^{\infty}{K_{p}(\cdot,x_{i})\tau_{i} }\colon\begin{subarray}{l}  (x_{i},\tau_{i}) \in \R^{d}\times \Lambda^{d-1}(\R^{d}\times  \R^{d}),\\ \norm{f}_{W}<\infty\end{subarray}\bigg\}
\end{equation}
where
\[
\norm{f}_{W}^{2}=\sum_{i,j=1}^{\infty}{K_{p}(x_{i},x_{j})\langle \tau_{i},\tau_{j}\rangle_{\Lambda^{d-1
    }(\R^{d}\times  \R^{d}) } }.
\]
The space $W$ contains forms varying only in $x$ (as $K_s=1$). Under this identification, the action of an embedded normal cycle in $N(S) \in W^{*}$ may be re-expressed as
\begin{align}\label{newaction}
        N(S)(\omega) = \int_{\mathcal{N}_{S}}{(\omega(x)\mid \tau^{\mathcal{N_{S}}}(x,\nu))\,d\mathcal{H}^{d-1}(x,\nu)},\quad \omega \in W.
\end{align}
This is the definition of $W$ and the action of $N(S) \in W^{*}$ that we use in \cref{mainsec}.
\end{remark}
\section{Decomposition and compression of normal cycles}\label{mainsec}
The compression algorithm of \cite{SparseNystromCurrVar} requires the target normal cycle in an explicit Dirac delta form. We now derive such a form for $N(\mathcal{T}) \in W^{*}$ when $\mathcal{T}$ is a discrete surface and $K_{s}=1$, matching the setting of~\cite{NCpaper}.
\begin{remark}
For discrete curves, the Dirac delta decomposition of the associated normal cycle is straightforward~\cite{NCThesis}: given discrete curves $C,S$ and $K_{s}=1$, one may represent
\begin{align*}
N(C)= \pi\sum_{i=1}^{n}{\delta_{x_{i}}A_{i} } \in W^{*},\quad N(S)=\pi\sum_{i=1}^{m} {\delta_{y_{i}}B_{i} } \in W^{*}
\end{align*}
where we use the notation of \cref{NCCURVEMET} for $x_{i},y_{i},A_{i},B_{i}$ and $W$ is identified with a space of vector fields induced by spatial kernel $K_{p}$.
\end{remark}


\subsection{Decomposition for discrete surfaces}
 Given a triangulation $\mathcal{T} = \cup_{i=1}^{n}{T_{i}} \subset \R^{3}$ (so $d=3$), a decomposition of $N(\mathcal{T})$ into distinct `planar', `cylindrical' and `spherical' components is derived in \cite{NCThesis} such that
\begin{align}\label{NC_decomp}
    N(\mathcal{T}) = N(\mathcal{T})^{\mathrm{\mathrm{pln}}} + N(\mathcal{T})^{\mathrm{\mathrm{cyl}}} + N(\mathcal{T})^{\mathrm{\mathrm{sph}}},
\end{align}   
and the three components are mutually orthogonal in $W^{*}$. The derivation given in \cite{NCThesis}  proceeds by recursively applying \cref{unionformula} to the union of triangles $\mathcal{T}=\cup_{i=1}^{n}{T_{i}}$ to obtain the decomposition. 
This decomposition leads to the following approximation.
\begin{lemma}\label[lemma]{mainlemma}
       Let $\mathcal{T}$ be a triangulated surface. The components of the associated discrete normal cycle embedding into $W^{*}$ (for $K_{s}=1$) may be expressed as
\begin{align*}
N(\mathcal{T})^{\mathrm{pln}}=0,\quad N(\mathcal{T})^{\mathrm{cyl}}
     = {2} \sum_{j=1}^{n_{e}}{\delta_{c_{j}}\gamma_{j} }+\sum_{j=1}^{n_e} R_j,\quad N(\mathcal{T})^{\mathrm{\mathrm{sph}}} =  \pi\sum_{x_{k} \in \partial \mathcal{T }}\delta_{x_{k}}\alpha_{k},
  \end{align*}
as elements of $W^{*}$, 
where the remainder
\begin{align*}
  R_j(\omega)=   \int_{e_j \times C}{ \omega_{x}(\tau_{x,u}) \,d\mathcal{H}^{2}(x,u) } -  \lvert e_j \rvert \int_{C}{\omega_{c_j}(\tau_{c_j,u})\,d\mathcal{H}^{1}(u) },
\end{align*}
and the coefficients
\begin{gather*}
    \alpha_{k} \coloneq \sum_{f_{i} \mathrm{bdry \ edge \ of}\ x_{k}}\sum_{j=1}^{3}{\frac{f_{ij}}{\lvert f_{i} \rvert} }\tilde{b}_{j} ,\quad \gamma_{j}\coloneq \begin{pmatrix}  f_{j}  \\ 0  \end{pmatrix} \wedge { \begin{pmatrix} 0\\   \underset{T\ \mathrm{triangle\,\, of}\ f_{j} }{\sum}{n_{T} } \end{pmatrix} },
\end{gather*} 
with $c_{i},x_{k},f_{i},n_{T,f_{i}}$ defined in the notation of \cref{NC_disc_met}, and $f_{ij}$ denoting the $j$th component of $f_i$, and
\begin{align*}
 \{\tilde{b}_{1},\tilde{b}_{2},\tilde{b}_{3}\}=\bigg\{\begin{pmatrix} 0 \\ e_2  \end{pmatrix}\wedge  \begin{pmatrix} 0\\ e_3 \end{pmatrix},\begin{pmatrix} 0 \\ e_1  \end{pmatrix}\wedge  \begin{pmatrix} 0\\ e_3 \end{pmatrix} ,\begin{pmatrix} 0 \\ e_1  \end{pmatrix}\wedge  \begin{pmatrix} 0\\ e_2 \end{pmatrix}   \bigg\}
\end{align*}
such that $e_{i}$ are standard basis elements of $\R^{3}$. 

\end{lemma}
The proof may be found in \cref{APPENDIXSECPROOF}, and relies on the computations of \cite{NCpaper,NCThesis}. 
The remainder $R_j(\omega)\coloneq \mathcal{C}(\omega) - \hat{\mathcal{C}}(\omega)$ 
associated to edge $e_j$ satisfies $|R_j(\omega)| = \mathcal{O}(|e_j|^2)$ 
for $\omega \in C^1$ by an analogous argument 
to~\cite[Proposition~1]{Kaltenmark2017AGF}, and vanishes as 
$|e_j| \to 0$. This regularity holds for all $\omega \in W$, 
as the Gaussian kernel ensures elements of $W$ are smooth. 
This is the same midpoint approximation used 
in~\cite{NCpaper} to derive \cref{thm:NCmetric}. In the remainder of the paper, we neglect the error term and write
\begin{equation}\label{NC_delta}
N(\mathcal T)=2\sum_{j=1}^{n_e}\delta_{c_j} \gamma_j + \pi\sum_{x_k\in\partial T}\delta_{x_k}\alpha_k.
\end{equation}



\subsubsection{Vector-valued embedding}
With the Dirac delta decomposition \cref{NC_delta}, we are almost in a position to apply the compression algorithm of \cite{SparseNystromCurrVar} to obtain a sparse approximation to $N(\mathcal{T}) \in W^{*}$. However, the compression algorithm assumes that $W$ is an RKHS of $\R^{k}$ valued functions for some $k\in \N$. To address this, we fix a standard basis of $\R^{6}$ and $\Lambda^{2}(\R^{3}\times \R^{3}) \simeq \Lambda^{2}(\R^{6})$ (of dimension $\binom{6}{2}= 15$), and use the resulting basis map to construct an isometric isomorphism of $W$ to an RKHS $V$ of $\R^{15}$-valued functions. 
This type of embedding arises in \cite{FidelityMetrics} for currents and varifolds, where $V$ is an RKHS of vector fields and real-valued functions.

Concretely, one begins by fixing a choice of isomorphism $B\colon \Lambda^{2}(\R^{6}) \to \R^{15}$ induced by the standard basis $\tilde{\mathcal{B}}=\{e_{i}\}_{i=1}^{6}$ of $\R^{6}$. In particular, given a wedge product $ a \wedge b \in \Lambda^{2}(\R^{6})$, we may expand in terms of $\tilde{\mathcal{B}}$ to obtain
\begin{align}
   a\wedge b = \bigg(\sum_{i=1}^{6}{a_{i}e_{i}}\bigg)\wedge \bigg(\sum_{j=1}^{6}{b_{j}e_{j}}\bigg) = \sum_{i,j=1}^{6}{a_{i}b_{j} (e_{i}\wedge e_{j})}\notag\\ = \sum_{1\le i<j\le 6}{(a_{i}b_{j} - a_{j}b_{i}) (e_{i}\wedge e_{j})} \eqcolon \sum_{k=1}^{15}{c_{k}(a,b)\beta_{k}},\label{wedging2}
\end{align}
using the bilinearity and antisymmetry properties of the wedge product. In \cref{wedging2}, we have labelled the corresponding orthonormal wedge basis
$
    \mathcal{B} \coloneq \{\beta_{k}\}_{k=1}^{15} \coloneq \{e_{i}\wedge e_{j}\}_{1\le i<j\le 6}
$
induced by $\tilde{\mathcal{B}}$. The coefficients $c_{k}(a,b)$ are the antisymmetric coefficients in \cref{wedging2} of the $k$th wedge basis element $\beta_k$. With respect to these fixed standard bases $\mathcal{B},\tilde{\mathcal{B}}$, one may define $B\colon \Lambda^{2}(\R^{6}) \to \R^{15}$ as the linear basis map defined by 
$
  B(\beta_{k})=\tilde{e}_{k},\quad  k\in \{1,\dots,15\},   
$
where $\tilde e_{k}$ denote standard basis vectors of $\R^{15}$ (distinct from the basis $e_i$ of $\R^6$). This mapping is extended linearly to general wedge products as
\begin{align*}
   B(a\wedge b) = \sum_{k=1}^{15}{c_{k}(a,b)\tilde e_{k}}  \in\R^{15}.
\end{align*}
 As $B$ is a linear map preserving orthonormal bases of the two finite-dimensional inner product spaces, it is an isometric isomorphism. 


Having fixed $B$ (with respect to choice of standard bases), we may isometrically identify the embedding space $W$ in \cref{newembeddingspace} of $\Lambda^{2}(\R^{3}\times \R^{3})$ valued mappings, with an RKHS $V$ of $\R^{15}$-valued functions, induced by the kernel
\begin{align*}
    K((x,u),(y,v) ) = K_{p}(x,y)\mathrm{I}_{{15}},\quad \forall (x,u),(y,v) \in \R^{3}\times \mathbb{S}^{2}.
\end{align*}
This identification is induced by an isometric isomorphism $F\colon W \to V$ with action
$
    F(K_{p}({x},\cdot)\tau) \coloneq K_{p}({x},\cdot)B(\tau),
$
which extends linearly to general elements of $W$. This in turn induces an isometric isomorphism of dual spaces $\tilde{F}\colon W^{*} \to V^{*}$, which allows one to embed $N(\mathcal{T})$ of \cref{NC_delta} into $V^{*}$ as
\begin{align}\label{realembedding}
    \tilde{F}(N(\mathcal{T})) = 2\sum_{i=1}^{n_{e}}{ \delta_{c_{i}} B(\gamma_{i})   } + {\pi}\sum_{x_{i} \in \partial\mathcal{T}}{\delta_{x_{i}} B(\alpha_{i})} \in V^{*}
\end{align}
with coefficient vectors in $\R^{15}$. Given the weights $\gamma_{i},\alpha_{i}$ associated to $N(\mathcal{T})$ in \cref{NC_delta}, each of the weights $B(\gamma_{i}), B(\alpha_{i})$ of $\tilde{F}(N(\mathcal{T}))$ can be efficiently computed in $\mathcal{O}(1)$. With the embedding \cref{realembedding}, one is now in a position to apply the compression algorithm of \cite{SparseNystromCurrVar} to approximate $\tilde{F}(N(\mathcal{T})) \in V^{*}$. Since $\tilde{F}$ is an isometric isomorphism, forming approximations to $\tilde{F}(N(\mathcal{T}))$  in $V^{*}$ will yield approximations of the same quality to $N(\mathcal{T}) \in W^{*}$, by applying $\tilde{F}^{-1}$.

\subsection{Compression of normal cycles}
The decomposition \cref{NC_delta} (and its curve analogue) writes the embedded normal cycle in the generic form
\begin{align*}
    \mu = \sum_{i=1}^{n}{\delta_{x_{i}}\alpha_{i}} \in V^{*},\quad (x_{i},\alpha_{i}) \in \R^{d}\times \R^{k},
\end{align*}
with $k=15$ for triangulated surfaces in $\R^{3}$, where $V$ is the RKHS of $\R^{k}$-valued functions induced by the product kernel $K_{p}\,\mathrm{Id}_{k}$. The compression algorithm of \cite{SparseNystromCurrVar} produces a sparse approximation $\hat{\mu}=\sum_{i=1}^{m}{\delta_{c_{i}}\beta_{i}}$ with $m\ll n$, as follows.
\begin{algorithm}[H]
\caption{Discrete functional compression with RLS sampling}
\begin{algorithmic}[1]
\State Fix domain $\R^{d}$, $m \ll n$, RKHS $V$, RKHS kernel function $K_{p}:\R^{d}\times \R^{d} \longrightarrow \R$ and target functional
\begin{align*}
    \mu = \sum_{i=1}^{n}{\delta_{x_{i}}\alpha_{i}} \in V^{*},\quad (x_{i},\alpha_{i}) \in \R^{d}\times \R^{k}.
\end{align*}
\State Sample $m$ distinct control points $\{c_{i}\}_{i=1}^{m} \subset \{x_{i}\}_{i=1}^{n}$ using approximate ridge leverage score (RLS) \cite{Musco,DAC} sampling.
\State Compute the vectors
\begin{align*}
y_{j} = \sum_{i=1}^{n}{K_{p}(c_{j},x_{i})\alpha_{i}} \in \R^{k},\quad j=1,\dots,m,
\end{align*}
and define $y=[y_{1},\dots,y_{m}]^{T} \in \R^{m\times k}$.
\State Form the control point approximation via orthogonal projection
\begin{align*}
    \hat{\mu} = \sum_{i=1}^{m}{\delta_{c_{i}}\beta_{i}} \in V^{*},\quad \beta = K_{CC}^{-1}y \in \R^{m\times k},\quad [K_{CC}]_{ij}=K_{p}(c_{i},c_{j})
\end{align*}
\end{algorithmic}\label{algorithm1}
\end{algorithm}
Control point sampling using approximate RLS sampling algorithms of \cite{Musco,DAC} (with $\mathcal{O}(m^{2}n)$ complexity) is fast in practice and yields strong theoretical bounds for decay of the compression error, as proven in \cite{SparseNystromCurrVar}. 
\begin{thm}[\cite{SparseNystromCurrVar}]\label{oldthm}
     Suppose we have a discrete target functional of the form
    \begin{align*}
        \mu_{S}=\sum_{i=1}^{n}{\delta_{x_{i}}\alpha_{i}}  \in V^{*},\quad (x_{i},\alpha_{i}) \in \R^{d}\times \R^{k}
    \end{align*}
with associated dual vector-valued function
\begin{align*}
   v^{\alpha}(x) = \sum_{i=1}^{n}{K_{p}(x,x_{i})\alpha_{i}} .
\end{align*}
Subsample $m$ distinct control points $\bm{c}=\{c_{i}\}_{i=1}^{m} \subset \{x_{i}\}_{i=1}^{n}$ and define the matrix of values 
    \begin{align*}
        Y_{C} = (v^{\alpha}(c_{1}),\dots,v^{\alpha}(c_{m}))^{T} \in \R^{m\times k},
    \end{align*}
    which is the evaluation of the dual vector-valued function on the control point locations. 
    Computing weights
    \begin{align*}
        \beta = [\beta_{1},\dots,\beta_{m}]^{T}= K_{CC}^{-1}Y_{C} \in \R^{m\times k},\quad [K_{CC}]_{ij}=K_{p}(c_{i},c_{j})
    \end{align*}
    yields an approximation
    \begin{align*}
       \hat{\mu}_{S}=   \sum_{i=1}^{m}{\delta_{c_{i}}\beta_{i} } 
    \end{align*} 
    that satisfies
       \begin{align}\label{tracebound}
      \norm{\mu_{S}-\hat{\mu}_{S} }_{V^{*}}^{2} \leq C\operatorname{Tr}(K_{XX}-Q_{XX}),\quad Q_{XX}\coloneq K_{XC}K_{CC}^{-1}K_{CX}.
    \end{align}

Furthermore, fixing $\delta\in (0,\frac{1}{32})$ and $S\in \N$, by subsampling control points using the recursive RLS scheme of \cite{Musco}, with probability $1-3\delta$ we have for a randomized sample size $m\approx  S\log(S/\delta)$ that
   \begin{align}\label{eigenbound}
        \norm{{\mu_{S} - \hat{\mu}_{S} }}_{V^{*}}^{2} \leq \frac{Cn}{S}\sum_{i=S+1}^{n}{\lambda_{i}(K_{XX})}.
    \end{align}
\end{thm}

That is, the randomised compression error decays at least as fast as the eigenvalues of the RKHS kernel matrix, evaluated on the original delta centres. For typical choices of the spatial kernel such as Gaussian RBF, this decay is exponential \cite{SparseNystromCurrVar} in $m$, and allows one to compress the target normal cycle to a very small level $m\ll n$, with negligible error in the RKHS norm. We demonstrate this in the following section.

\section{Numerical experiments}\label{numexp}
As in our work \cite{SparseNystromCurrVar} for currents and varifolds, we now demonstrate the effectiveness of \cref{algorithm1} for compression of large-scale embedded normal cycles when $K_{s}=1$, on surfaces from modern geometry processing datasets.\footnote{ Python code is available at \url{https://github.com/tonyshardlow/GeometricMeasureCompression}.} We demonstrate properties such as rapid approximation error decay and effectiveness in downstream applications to nonlinear LDDMM registration. All experiments in this section are performed on a Tesla T4 GPU with 16GB of RAM. We also use the PyKeops library \cite{Keops} to speed up all kernel reduction operations.
\subsection{Compression error decay}

\begin{figure}[H]
\begin{subfigure}{1.\textwidth}  
    \includegraphics[width=.4\textwidth]{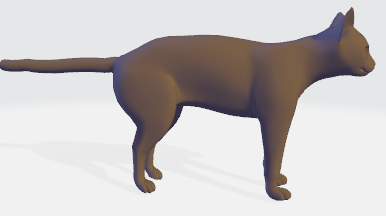}
    \hfill
    \includegraphics[width=.2\textwidth]{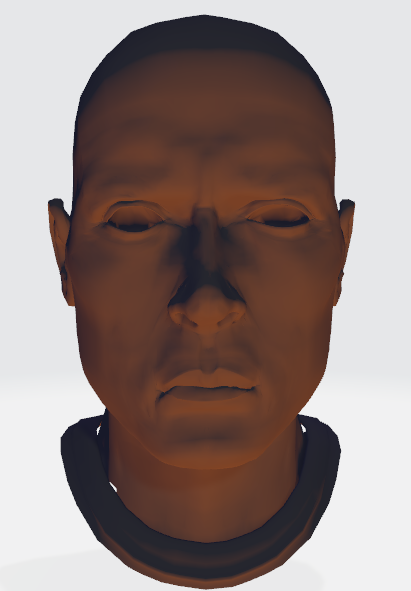}
    \hfill
    \includegraphics[width=.24\textwidth]{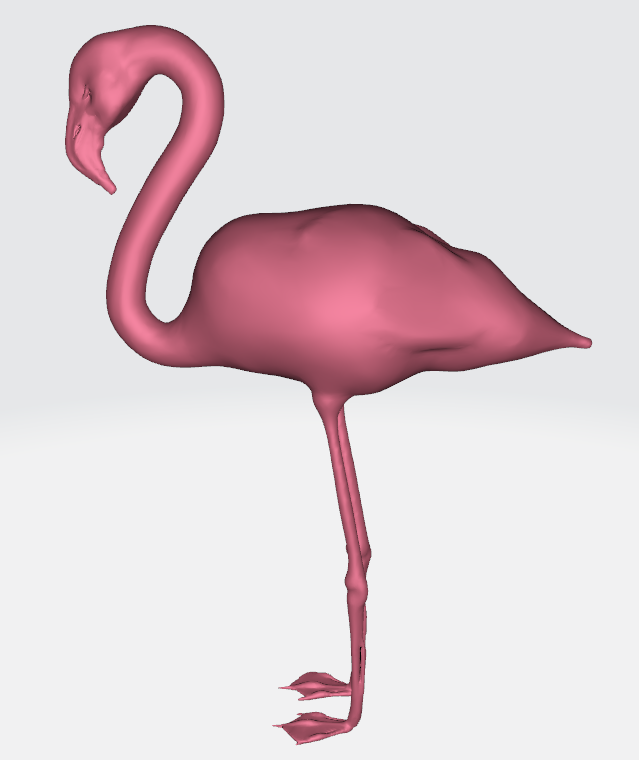}
\end{subfigure}
   \caption{ {Left}: Cat (14,410 triangles). {Middle}: Head (31,620 triangles). {Right}: Flamingo (52,895 triangles).}\label{firstdatfig}
\end{figure}

We begin by illustrating the decay of the normal cycle compression error of \cref{algorithm1} in the RKHS dual norm as a function of $m$, the compression size. We use the surface data shown in \cref{firstdatfig} for our experiments, originally developed by the authors of \cite{ANIMdataset} and used in our prior work \cite{SparseNystromCurrVar} for currents and varifolds.
The data are centred and scaled, so that the flamingo surface lies in a box of size $1.3\times 3.3 \times 7.1$, the head surface lies in a box of size $3.8\times 5.3 \times 4.0$, and the cat surface lies in a box of size $1.5\times 5.3 \times 3.8$. 

\begin{figure}[H]
\centering  
\begin{subfigure}{\textwidth}

\includegraphics[width=.34\textwidth] {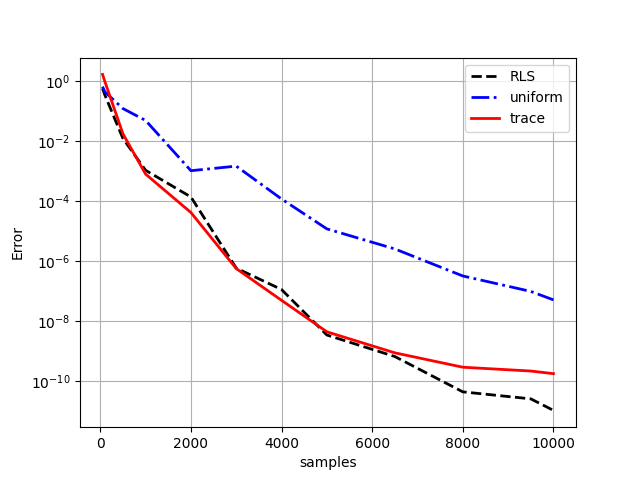}\includegraphics[width=.34\textwidth]{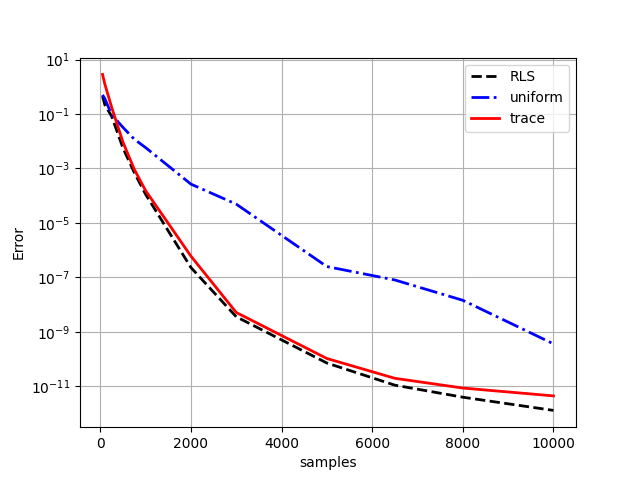}\includegraphics[width=.34\textwidth]{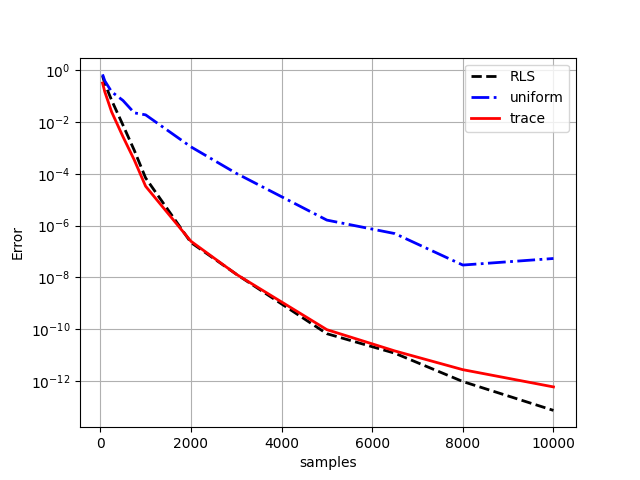}
\end{subfigure}
   \caption{Numerical curves comparing RKHS (relative) error decay (black) of Algorithm~1 for compression of normal cycles, to theoretical trace bound (red) and uniformly sampled compression (blue), on cat (left), head (centre) and flamingo (right) surfaces.  }
    \label{fig:curves}
\end{figure}

For each test surface, we run Algorithm~1 and plot the true (relative) square error in RKHS norm of the compression as a function of $m$. We also plot for comparison the error curve for uniform sampling, where control points are uniformly subsampled without replacement and used as delta centres for orthogonal projection. Finally, the trace bound from \cref{oldthm} is also plotted with rescaled constants, for  comparison of the error decay rate. For the kernel of the embedding RKHS, we make a standard choice of Gaussian RBF kernel, $K_{p}(x,y)=\exp(-\frac{\norm{x-y}^{2}}{2\sigma^{2}})$, with scale parameters $\sigma\in \{0.3,0.5,0.25\}$ respectively, for the cat, head and flamingo test cases. The scales are chosen to be roughly the typical scale of geometric variations of interest for each shape. The resulting plots are shown in \cref{fig:curves}.

We observe in all cases that the compression error decays  quickly, allowing one to achieve high-quality approximations for $m \ll n$ across all examples. The theoretical trace-bound decay rates closely match those of the compression error decay of \cref{algorithm1}, which consistently outperforms uniform sampling. We also note that the compression time is short; even for the large flamingo test case, the compression to $1000$ delta centres ($<0.1\%$ error and around $99\%$ compression ratio) takes less than $1$ second. In practice, downstream tasks such as shape registration only require accuracy comparable to the intrinsic uncertainty of the data~\cite{DURRLEMAN}, which the curves of \cref{fig:curves} achieve at $m\ll n$.

\subsection{Matching quality}\label{matchqualitysec}
We now demonstrate the effectiveness of normal cycles compression for nonlinear registration of surfaces; a similar demonstration was given in our work \cite{SparseNystromCurrVar} for varifolds.

  As a baseline, we consider a Hamiltonian control-point LDDMM \cite{Sparseadaptiveparam} registration problem of a template surface $T  \subset \R^{d}$ to target $S \subset \R^{d}$ using the normal cycles metric. This problem involves a set of fixed spatial deformation control points $\bm{c}^{s}=(c_{i}^{s})_{i=1}^{P} \in \R^{Pd}$ subsampled from the $n_{v}$ vertices of the template, a smooth positive definite scalar spatial kernel $k$, and minimises an objective $E\colon \R^{Pd}  \to \R_{\geq 0}$, defined by the following system
\begin{align}\label{controlobjectiveintro}
E(\bm{\alpha}) \coloneq \frac{1}{2}\norm{v^{\bm{\alpha} }(0)}_{V_{k}}^{2} + \frac{\lambda}{2}\norm{\mu_{\varphi^{v^{\bm{\alpha} }}_{01}(T) } - \mu_{S} }_{V^{*}}^2,\quad \bm{\alpha} =(\alpha_{1},\dots,\alpha_{P} ) \in \R^{Pd},
\end{align}
where
\begin{align}\label{controlpointfieldintro}
v^{\bm{\alpha} }(t,x) \coloneq \sum_{i=1}^{{P}}{k(x,c_{i}^{s}(t))\alpha_{i}(t)},
\end{align}
and, for $i=1,\dots,P$,
\begin{align}\label{controlpointHamiltonianintro}
\begin{cases}\displaystyle
 \partial_{t}c_{i}^{s}(t) = \sum_{j=1}^{P}{k(c_{i}^{s}(t),c_{j}^{s}(t))\alpha_{j}(t)},\qquad c_{i}^{s}(0) = c_{i}^{s},\\ \displaystyle
 \partial_{t}\alpha_{i}(t) = -\sum_{j=1}^{P}{\nabla_{1}k(c_{i}^{s}(t),c_{j}^{s}(t))\alpha_{i}(t)^{T}\alpha_{j}(t)},\qquad \alpha_{i}(0) = \alpha_{i}.
\end{cases}
\end{align}
 The initial momenta $\bm{\alpha}\in \R^{Pd}$ attached to the $P$ deformation control points induce a Hamiltonian motion \cref{controlpointHamiltonianintro}, a time-dependent vector field $v^{\alpha} \in L^{2}([0,1];V_{k}) $ (where $V_{k}$ is the RKHS of vector fields induced by $k$), and a flow $\varphi_{01}^{v^{\alpha}}$ deforming the template. It is common practice \cite{GlaunesCurrents} to approximate $\varphi(T)$ by computing the pushforward of the $n_{v}$ vertices of $T$ and forming a new triangulation from the deformed vertices with the \tit{same} mesh connectivity as the original template. This approximation has negligible error for finely sampled shapes. The initial momenta may be optimised with respect to  \cref{controlobjectiveintro}, so that $\mu_{\varphi^{v^{\bm{\alpha} }}_{01}(T) } \approx \mu_{S}$  in the dual metric on the normal cycle embedding.
 
 Computing the objective and its gradients splits into two parts: flow computation and metric computation. At each step of the ODE solver for the flow, one computes a reduction at cost $\mathcal{O}(Pn_{v} )$, which may be implemented efficiently using KeOps, and is fast as one takes $P\ll n_{v}$ in the control point parametrisation of deformations. The number of vertices $n_{v}$ is also much fewer than the number of edges $n_{e}$, and one usually only requires simple first or second-order solvers with few evaluations for the flow. Both in terms of complexity and memory, the metric (and its gradient) computation at $\mathcal{O}(n_{e}^{2})$ becomes the dominating cost; for very large $n_{e}$, this computation can become very slow even with KeOps.

 We compare this to a modified \tit{compressed} version of the above matching algorithm, which minimises the following objective
    \begin{align}\label{modifiedobjective}
E(\bm{\alpha}) \coloneq \frac{1}{2}\norm{v^{\bm{\alpha} }(0)}_{V_{k}}^{2} + \frac{\lambda}{2} \norm{\hat{\mu}_{S} - \hat{\mu}_{\varphi^{v^{\bm{\alpha} }}_{01}(T)} }_{V^{*}}^{2} ,\quad \bm{\alpha} =(\alpha_{1},\dots,\alpha_{P} ) \in \R^{Pd},
\end{align}

where the deformation control points and Hamiltonian dynamics \cref{controlpointHamiltonianintro} are unchanged between objectives \cref{controlobjectiveintro} and \cref{modifiedobjective}. The cost of flow (and flow gradient) computation remains the same as the deformation mechanics are unchanged. The main difference here is the compressed target normal cycle $\hat{\mu}_{S}$ and compressed deformed template normal cycle $\hat{\mu}_{\varphi^{v^{\bm{\alpha} }}_{01}(T)}$ in the discrepancy term. The target normal cycle is compressed once, offline (pre-optimisation), and the deformed template normal cycle is compressed at each iteration. We assume a fixed level of compression $m\ll n_{e}$ has been chosen. For the deformed template, instead of RLS sampling for compression at each iteration of the registration algorithm,  we sample the template normal cycle \textit{once pre-optimisation}, and project the deformed template at each iteration onto the transported RLS samples.

The cost of performing this projection of the deformed template per-iteration is $\mathcal{O}(m^{3} + mn_{e})$, which comes from the linear solve and reduction operation during compression. The reduction operation of cost $\mathcal{O}(mn_{e})$ may be implemented efficiently via KeOps and is rapid. The linear systems to be solved at cost $m^{3}$ are of size $m\ll n_{e}$, which in practice are small enough to fit comfortably on GPU memory, allowing one to parallelise and significantly accelerate the resulting linear solve. Therefore, compressing the deformed template normal cycle and computing the metric (and its gradients) using the compressed normal cycles is much faster than computation of the full metric and gradient. One may also achieve further speed-up by exploiting the structure of the systems to be solved (e.g., sparsity, symmetry) and using conjugate gradient type schemes, which have been recently combined with the efficient KeOps reduction framework \cite{Keops}.

We show that one may obtain high-quality registration results using the compressed matching algorithm above, even when only $1$--$2\%$ of the Dirac delta centres are used. Furthermore, we shall observe how one may obtain significant computational savings using the compressed version of the registration algorithm. We demonstrate this on two  densely sampled shapes, taken from a modern geometry processing database of surfaces \cite{Thingi10k}. The first is the Queen test surface with  {$229,776$} triangles, and the second is the PumpkinHead test surface with  {$394,510$} triangles. In the full normal cycle representation, the Queen test surface has $344,664$ Dirac delta centres, and the PumpkinHead surface has $591,765$ Dirac delta centres. The PumpkinHead surface lies in an origin-centred box of size $3.5\times 3.8 \times 3.9$, and the Queen test surface in a box of size $3.1 \times 4.0 \times 4.9$. 

In both test cases, we compare the registration quality and runtime of the registration algorithm for the uncompressed and compressed matching algorithms. We choose a spherical template, with the same triangle resolution as the target for both experiments, lying in a box centred at the origin of size $3.5\times 3.5 \times 3.5$. For the normal cycles representation, in both cases we choose $K_{p}$ as the Gaussian RBF kernel with spatial length-scale $\sigma_{p} = 0.2$. For the spatial kernel parametrising deformation vector fields in \cref{controlpointfieldintro}, we fix a sum of $4$ Gaussian kernels of decreasing length-scales $\sigma_{i} \in \{1.0,0.5,0.2,0.1\}$. For all experiments, the Hamiltonian systems \cref{controlpointfieldintro} are solved via a forward Euler scheme with $10$ time-steps. The kernel reduction and gradient computations for flows and both metric computations are performed using Keops and automatic differentiation. Optimisation of the LDDMM registration objective is performed via an LBFGS routine run for $500$ iterations, with strong-Wolfe line search used for setting learning rate. We set the number of spatial deformation control points $P$ and measure compression control points $m$ to be the same in all subsequent experiments. The final registration quality is measured in the Hausdorff metric over sets given as
\begin{align*}
    d_{H}(A,B) \coloneq \mathrm{max}\bigg(  \underset{a \in A}{\sup }{\ d(a,B)},  \underset{b \in B}{\sup }{\ d(b,A)}\bigg),\quad d(x,A) \coloneq \underset{a \in A}{\inf}{\norm{x-a}_{2}},
\end{align*}
which is independent of the optimisation objectives.

\begin{figure}[H]
 \hspace{50pt}\begin{subfigure}{0.4\textwidth}
     \includegraphics[width=\textwidth]{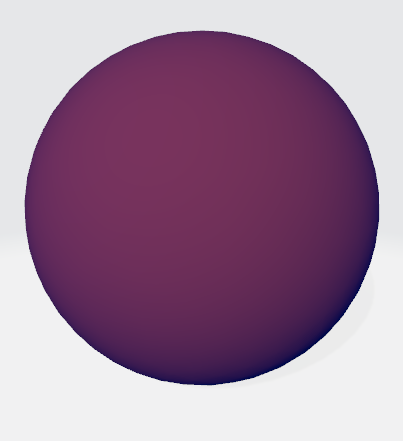}
 \end{subfigure}
 \hspace{20pt}
 \begin{subfigure}{0.35\textwidth}
     \includegraphics[width=\textwidth]{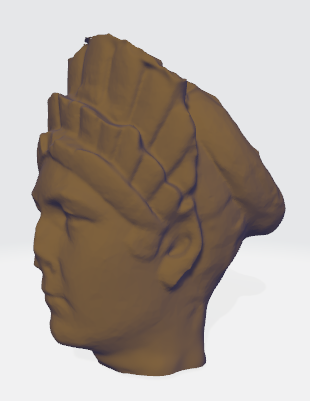}
 \end{subfigure}
 
 \medskip
 \hspace{70pt}\begin{subfigure}{0.35\textwidth}
     \includegraphics[width=\textwidth]{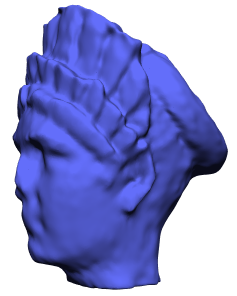}
 \end{subfigure}
 \hspace{30pt}
 \begin{subfigure}{0.35\textwidth}
     \includegraphics[width=\textwidth]{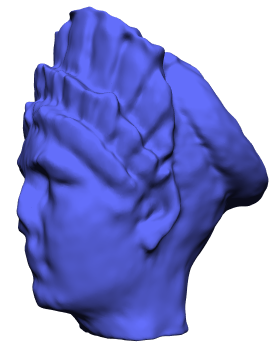}
 \end{subfigure}
 \caption{ 
   {Top left}: spherical template.  {Top right}: target mesh.   {Bottom left}: matching with full normal cycles taking  $2$ hours and $42$ minutes with $d_{H} = 0.0442$.  {Bottom right}: matching with $97\%$ compression of template and target taking only $17$ minutes with $d_{H} = 0.0298$.} \label{Queencompression}
\end{figure}

We begin with the Queen test surface, where the target and deformed template normal cycles are compressed down to $10,000$ Dirac deltas each, from the original $344,664$ Dirac delta centres. This is a compression ratio of over {$97\%$}. Data and results from the compressed and uncompressed matching are shown in \cref{Queencompression}.
We observe in \cref{Queencompression} that the matching quality is almost identical between compressed and uncompressed cases. This is numerically evidenced by the Hausdorff metric registration error for full ($d_{H} = 0.0442$) and compressed ($d_{H} = 0.0298$) normal cycle registration, with the compressed matching yielding a slightly better quality of registration. Importantly, the compressed matching algorithm yields a significant speed-up of $9$--$10$ times over the uncompressed version, reducing overall matching time from $9720$s to $1020$s.

\begin{figure}[H]
 \hspace{50pt}\begin{subfigure}{0.35\textwidth}
     \includegraphics[width=\textwidth]{matching/sphere.png}
 \end{subfigure}
 \hspace{20pt}
 \begin{subfigure}{0.35\textwidth}
     \includegraphics[width=\textwidth]{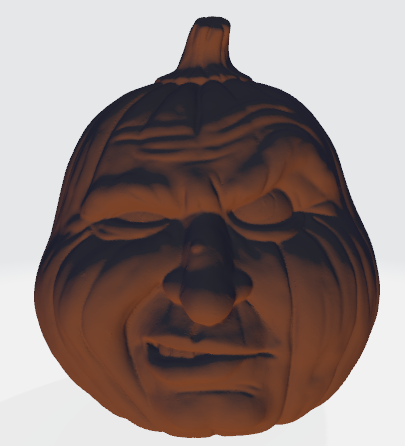}
 \end{subfigure}
 
 \medskip
 \hspace{15pt}\begin{subfigure}{0.52\textwidth}
     \includegraphics[width=\textwidth]{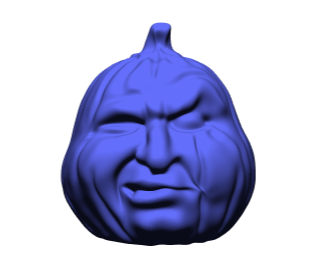}
 \end{subfigure}
 \hspace{-10pt}
 \begin{subfigure}{0.39\textwidth}
     \includegraphics[width=\textwidth]{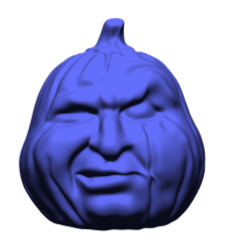}
 \end{subfigure}
 \caption{  {Top left}: spherical template.  {Top right}: target mesh.  {Bottom left}: an example matching with full normal cycles, from spherical template to target without compression, taking $5$ hours and $37$ minutes with $d_{H} = 0.1621$.  {Bottom right}: result of compressed matching problem with $99\%$ compression ratio, taking only $17$ minutes with $d_{H} = 0.1662$. } \label{Pumpkinfig}
\end{figure}

We now move on to the second test case on the Pumpkinhead test surface with {$394,510$} triangles, and $591,765$ Dirac centres in the normal cycle. The target and deformed template normal cycle are compressed down to $5000$ Dirac delta centres each from the original discrete normal cycles, yielding a compression ratio of over {$99\%$}. The data and matching results are shown in \cref{Pumpkinfig}.
 Once again, we obtain significant runtime savings, with a factor of $20$ times speed-up over the uncompressed matching problem, reducing overall matching time from $20,220$s to $1020$s. We also observe in \cref{Pumpkinfig} that a comparable visual quality of matching is achieved even with the severe level of compression. This is also reflected in the negligible difference in Hausdorff metric scores between the full ($d_{H} = 0.1621$) and compressed ($d_{H} = 0.1662$) matching problems.    
\subsection{Comparison with target-mesh decimation}\label{decimsec}
A widely used baseline in LDDMM and elastic shape-matching pipelines is to
\emph{decimate} the target mesh prior to registration, and then construct
the normal cycle of the decimated mesh, rather than compress the full
normal cycle directly. We compare these two strategies on the Queen and
PumpkinHead targets of \cref{matchqualitysec} at an equal budget of
$m=2{,}000$ Dirac centres each --- a more aggressive compression than the
$10{,}000$ and $5{,}000$ budgets used in \cref{Queencompression} and
\cref{Pumpkinfig}, chosen here so that any quality gap between the two
strategies is sharply visible. For the compression, we apply \cref{algorithm1} with ridge-leverage-score
sampling at $\sigma_{p}=0.2$, matching the kernel scale used in
\cref{matchqualitysec}. For the decimation, we run quadric edge-collapse
decimation~\cite{Garland1997QEM} as implemented in
\texttt{pymeshlab}~\cite{Pymeshlab}, with a face budget of $1{,}333$; for
these closed-manifold targets this yields $\approx 2{,}000$ edge-midpoint
Dirac centres, matching the compression budget. The resulting sparse normal
cycles $\hat{\mu}$ are then assessed directly against the full-resolution
target normal cycle $\mu_{S}$ along three axes: the relative
dual-norm discrepancy
$\|\hat{\mu}-\mu_{S}\|_{V^{*}}^{2}/\|\mu_{S}\|_{V^{*}}^{2}$, which is the
quantity \cref{controlobjectiveintro} minimises during LDDMM; the symmetric
mean Chamfer distance~\cite{Chamferref}, which captures average surface
coverage; and the symmetric Hausdorff distance $d_{H}$, which captures
worst-case point-wise deviation. Numerical results are given in
\cref{tab:decim}.

\begin{table}[h]
\caption{Comparison of sparse normal-cycle representations of the Queen and
PumpkinHead targets at equal budget $m=2{,}000$. Compression via
\cref{algorithm1} (Comp.) is roughly one and a half orders of magnitude more
accurate in the dual-norm discrepancy that the LDDMM objective
\cref{controlobjectiveintro} optimises, and is also slightly more accurate
in mean Chamfer distance. Quadric edge-collapse decimation (Decim.) achieves
a smaller worst-case $d_{H}$ on shapes with substantial low-curvature regions. Best of each pair in bold.}\label{tab:decim}
\begin{tabular}{@{}l|cc|cc|cc@{}}
\toprule
& \multicolumn{2}{c|}{Rel.\ NC discrepancy} & \multicolumn{2}{c|}{Chamfer} & \multicolumn{2}{c}{$d_{H}$} \\
Target & Comp. & Decim. & Comp. & Decim. & Comp. & Decim. \\
\midrule
Queen        & $\mathbf{0.73\%}$ & $9.08\%$  & $\mathbf{0.0715}$ & $0.0771$ & $0.490$ & $\mathbf{0.476}$ \\
PumpkinHead  & $\mathbf{0.38\%}$ & $13.41\%$ & $\mathbf{0.0643}$ & $0.0677$ & $0.940$ & $\mathbf{0.523}$ \\
\bottomrule
\end{tabular}
\end{table}
 
The dual-norm discrepancy --- the loss actually minimised by
\cref{controlobjectiveintro} --- is between $12$ and $35$ times smaller for
compression than for decimation, with \cref{algorithm1} achieving 
sub-percent error against the full target normal cycle on both shapes at
this aggressive budget. Mean Chamfer error is also marginally smaller for
compression, indicating that ridge-leverage-score samples, despite their
non-uniform density, cover the target surface as faithfully as the vertices
retained by quadric decimation. The Hausdorff distance $d_{H}$ is the only
metric on which decimation outperforms compression: on the Queen the two
methods are essentially tied, but on the PumpkinHead --- whose surface
contains large low-curvature regions on the back and flanks ---
ridge-leverage-score sampling concentrates Dirac centres on the high-curvature
features at the front, leaving a flat patch under-represented, which the
 worst-case metric $d_{H}$ then magnifies.
The error-decay curves of
\cref{fig:curves} are consistent with this picture: the blue uniform-subsampling
curves --- the natural unweighted analogue of geometric decimation in the
normal-cycle representation --- decay substantially more slowly than
\cref{algorithm1} on all three test shapes.

\subsection{Downstream applications and scaling to complex shapes}\label{downstreamsec}
The compression scheme of \cref{algorithm1} has subsequently been
deployed within a learning-based diffeomorphic shape-matching pipeline
in our companion work \cite{CurvatureEnthusiasm}, providing further
empirical evidence of its practical utility on shapes with regions of
high curvature. There, compressed normal cycles drive the matching
loss for articulated meshes drawn from the TOPKIDS, DFAUST, MANO, SMAL and
TOSCA datasets, including human, hand, animal and child body
scans whose curvature profile is comparable to (and in places more
demanding than) the cat, flamingo and Stanford Bunny shapes
of~\cite{NCpaper, SparseNystromCurrVar}. In particular, Figure~13 of
\cite{CurvatureEnthusiasm} matches a $12$k-vertex source to a $120$k-vertex
target with normal cycles compressed to $2.5$k, $5$k and $10$k Dirac
centres, with Chamfer-error degradation that is visually
indistinguishable from the uncompressed match.

A detailed runtime breakdown for the compressed normal-cycle
representation is also provided in~\cite{CurvatureEnthusiasm}. Figure~15
of that work separates per-iteration metric/gradient cost from
total wall-clock time, and Figure~16 (also reported as raw timings in
the supplementary table) shows that, on DFAUST samples upsampled to
between $10$k and $500$k vertices via Instant Field-Aligned
Meshes~\cite{Jakob2015IFA}, wall-clock time for compressed normal
cycles ($5$k or $10$k Dirac centres) remains essentially \textit{flat}
in the target resolution, whereas the uncompressed cost grows
super-linearly. This is consistent with the analysis of complexity in
\cref{matchqualitysec}: once the target normal cycle has been
compressed (a one-off cost that we measured in \cref{fig:curves} to be
under one second for the largest test case), the dominant
$\mathcal{O}(m\,n_{e})$ reduction is independent of the original mesh
size. The earlier Varifold version of the same scaling experiment
appears as Figure~1 of~\cite{ARCFlow}, where the same qualitative
behaviour was first reported.

A complementary picture emerges for memory. The dense kernel block
$K_{CC} \in \R^{m\times m}$ required for metric and gradient
computation inside the registration loop of \cref{matchqualitysec} has
$m^{2}=10^{8}$ entries at the Queen budget of $m=10{,}000$ and
$2.5\times 10^{7}$ entries at the PumpkinHead budget of $m=5{,}000$,
both fitting comfortably in $16$~GB of GPU memory at single precision.
The uncompressed equivalents are of order $n_{e}^{2}\approx 1.2\times
10^{11}$ and $3.5\times 10^{11}$ entries respectively---three to four
orders of magnitude larger, well past the capacity of a single GPU, and
only tractable at all via the lazy reduction strategy of KeOps. The
compressed representation $\{(c_{i},\beta_{i})\}_{i=1}^{m}$ itself
occupies only $\mathcal{O}(m\,k)$ memory with $k=15$, i.e.\ on the
order of megabytes for the experiments considered here.

The ridge-leverage-score sampling step at the heart of
\cref{algorithm1} has a memory profile of its own which becomes
relevant on very large targets. The recursive RLS scheme of
\cite{Musco}, which we use in \cref{numexp}, has peak memory of order
$\mathcal{O}(m\,n)$ in its leverage-score-estimation step, where $n$
is the size of the discrete normal cycle being compressed; this is
comfortable for $n\lesssim 10^{6}$ on single-GPU hardware (covering
all the targets considered in this paper), but is the binding
constraint as $n$ grows further. The divide-and-conquer RLS scheme of
\cite{DAC} addresses this by partitioning the input into $\sqrt{n}$
disjoint batches of size $\sqrt{n}$, computing leverage scores
\textit{per batch} (so that only one batch is resident in memory at a
time), and merging the results; the peak working memory is therefore
$\mathcal{O}(m\sqrt{n})$, at the cost of a modest constant overhead in
the sampling complexity.
\cref{algorithm1} is compatible with either choice; the user may
select between them based on the scale of the target.

\section{Conclusion}
In this work, we have extended our algorithm in \cite{SparseNystromCurrVar} to compression of normal cycles, allowing us to scale up the normal cycles metric to shapes with large-scale resolution of the order $10^{5}$--$10^{6}$, while keeping the computational cost at a fraction of the full representation. Furthermore, our compression comes with fast error decay and guarantees of small RKHS distance to the normal cycle of the true underlying shape. We have demonstrated the effectiveness of the compression method on real-world shape data in downstream tasks such as nonlinear LDDMM registration. Indeed, our experiments show dramatically reduced run times, with negligible difference in registration quality. We leave as future work the extension of the compression algorithm to linear and Gaussian normal kernel cases, which can help increase the sensitivity of the underlying shape representations to higher-order curvature information.

\backmatter

\bmhead{Code availability}
The Python implementation accompanying this paper is openly available at \url{https://github.com/tonyshardlow/GeometricMeasureCompression}.

\bmhead{Data availability}
All shape data used in \cref{numexp} are taken from publicly available sources cited in the text: \cite{ANIMdataset} for the cat, head and flamingo meshes, and \cite{Thingi10k} for the Queen and PumpkinHead meshes.

\appendix

\section{Background}\label{Appendix1}
We briefly introduce  definitions required for defining currents and normal cycles in \cref{bg}.
    \subsection{Wedge product and differential forms}
Given two vectors $v,w \in \R^{d}$, the wedge  or exterior product is an operation on $(v,w)$ that produces an object $v\wedge w$ in a way that naturally generalises the properties of the cross product in $\R^{3}$, both algebraically and geometrically. In particular, $v\wedge w$ contains information about the oriented parallelogram spanned by $v,w$, with modulus equal to the area of the associated parallelogram. This product is characterised by the following properties for all $\alpha \in \R$, $v,w,z \in \R^d$:
\begin{description}\label{items}
    \item [$\wedge$ is bilinear:] $(v+\alpha w) \wedge z = v \wedge z + \alpha(w \wedge z)$,
    \item[$\wedge$ is alternating:] $v \wedge v = 0$, and
    \item[$\wedge$ is antisymmetric:] $v \wedge w = -w  \wedge v$.
\end{description}
 The space in which $v\wedge w $ lives is denoted $\Lambda^{2}\R^{d}$, the $2$nd exterior power of $\R^{d}$, which consists of formal linear combinations of elements of the form $v\wedge w$, giving it a vector space structure. One can also define analogously $\Lambda^{3}\R^{d}$, the third exterior power of $\R^{d}$ consisting of formal linear combinations of elements $v\wedge w\wedge z$, representing an oriented parallelepiped, spanned by the three vectors. This process can be repeated to obtain the $k$th exterior power of $\R^{d}$, denoted $\Lambda^{k}\R^{d}$ for all $k<d$. 

We now give a formal description of the exterior product spaces following \cite{Sternberg,FrankMorgan,NCThesis}. We begin by defining the space of alternating multilinear forms on $\R^{d}$.



\begin{definition}
    An alternating multilinear form of degree $k$ is a mapping 
    \begin{align*}
    f: \underbrace{\R^{d}\times \dots \times \R^{d}}_{k \ \mathrm{times}}\longrightarrow \R
    \end{align*}
    with the properties that $f$ is linear in each entry and  $f(v_{1},\dots,v_{k})=0$, whenever $v_{1},\dots,v_{k}$ are linearly dependent.
    Note that these two properties automatically imply 
    \begin{align*}
        f(v_{\sigma(1)},\dots,v_{\sigma(k)}) = \mathrm{sgn}(\sigma)f(v_{1},\dots,v_{k})
    \end{align*}
    for all permutations $\sigma$ of $\{1,\dots,k\}$. We denote the space of alternating $k$-multilinear forms as $A^{k}\R^{d}$, which is a finite-dimensional vector space under pointwise addition and multiplication. 
\end{definition}
One may subsequently define the $k$th exterior product space $\Lambda^{k}\R^{d}$ as follows.
\begin{definition}
For $0\leq k\leq d$, one defines the $k$th exterior product of $\R^{d}$ as $\Lambda^{k}\R^{d}\coloneq (A^{k}\R^{d})^{*}$, the dual of the space of alternating $k$-multilinear forms. One further defines the wedge product of vectors $v_{1},\dots,v_{k} \in \R^{d}$ as the element denoted $v_{1}\wedge \dots \wedge v_{k} \in \Lambda^{k}\R^{d}$ with action
\begin{align}\label{wedgedef}
    (v_{1}\wedge \dots \wedge v_{k})(f) \coloneq f(v_{1},\dots,v_{k}),\quad  f \in A^{k}\R^{d}.
\end{align}
It can be shown that products of this form span $\Lambda^{k}\R^{d}$, so that the dual space is made up of linear combinations of such elements. In particular, one has $\mathrm{dim}(\Lambda^{k}\R^{d}) =  \binom{d}{k}$ as a basis of this space is given by $\mathcal{B}=\{e_{i_{1}}\wedge\dots \wedge e_{i_{k}}\}_{1\le i_{1}<\cdots<i_{k}\le d}$. Furthermore, it can be verified that the wedge product of $k$ vectors as defined in \cref{wedgedef} is linear in each element, and vanishes when $v_{1},\dots,v_{k}$ are linearly dependent, thus satisfying the characterising properties.
\end{definition}


One may define an inner product on $\Lambda^{k}\R^{d}$ by defining it on simple wedges of $k$ vectors as
\begin{align}\label{exteriorinnerproddef}
    \langle v_{1}\wedge \dots \wedge v_{k},w_{1}\wedge \dots \wedge w_{k}\rangle_{(\Lambda^{k}\R^{d}) } \coloneq \det((\langle v_{i} ,w_{j}\rangle)_{1\leq i,j\leq k})
\end{align}
and extending linearly. One may verify that this inner product is a simple Euclidean inner product on components with respect to the basis $\mathcal{B}$. In \cref{bg,mainsec,APPENDIXSECPROOF}, we often use the Riesz representation theorem on the resulting inner product space to isometrically identify $(\Lambda^{k}\R^{d})^{*}\simeq (\Lambda^{k}\R^{d})$.


These concepts may be used to define differential forms. For the purposes of defining the current associated to an $m$-dimensional submanifold, one usually considers compactly supported differential $m$-forms, defined as
$
  \Omega_{0}^{m}(\R^{d}) \coloneq C_{0}(\R^{d},(\Lambda^{m}\R^{d})^{*} )
$
equipped with the infinity norm 
\begin{align*}
    \norm{\omega}_{\infty} = \sup_{x\in \R^{d}}{\lvert \omega(x)\rvert}_{(\Lambda^{m}\R^{d})^{*}}. 
\end{align*}
A differential $m$ form on $\R^{d}$ is a natural object that one can integrate over oriented $m$-dimensional manifolds, in a coordinate-independent manner.

\subsection{Hausdorff measure and rectifiable sets}
The general setting in which one defines currents is that of rectifiable sets with respect to Hausdorff measure. We briefly define these concepts here.
\subsubsection{Hausdorff measure}
The $m$-dimensional Hausdorff measure on $\R^{d}$ is defined as follows. First, one defines
\begin{align*}
    \mathcal{H}_{\delta}^{m}(S) \coloneq \inf\bigg\{ \sum_{i=1}^{ \infty}{\mathrm{diam}(U_{i})^{m}}\colon U_{i}\ \mathrm{open},\  S \subset \cup_{i=1}^{\infty}{U_{i}},\quad \mathrm{diam}(U_{i}) < \delta\bigg\},
\end{align*}
which is the effective size of the smallest cover by open sets with diameter bounded by $\delta$. One takes the limit as $\delta \to 0$, to obtain the outer measure
\begin{align*}
    \mathcal{H}^{m}(S) \coloneq \lim_{\delta\to 0}{\mathcal{H}^{m}_{\delta}}(S),
\end{align*}
which can be made into a measure by restricting to measurable sets with respect to $\mathcal{H}^{m}$. The resulting measure is $m$-dimensional Hausdorff measure on $\R^{d}$. When $m=d$, this agrees with $\lambda^{d}$, the Lebesgue measure in $\R^{d}$. When $m<d$, this is a natural measure to use for integration and measuring `volume' on $m$-dimensional subsets $S$ embedded in $\R^{d}$, on which the Lebesgue measure would be zero. The Hausdorff measure may be used to measure sub-dimensional volume for very general subsets, and $\mathcal{H}^{m}(S)$ agrees with the classical/intuitive notion of volume for submanifolds $S$ of dimension $m$ embedded in $\R^{d}$. One may also define the Hausdorff dimension of a set as 
\begin{align*}
    \mathrm{dim}_{\mathcal{H}} = \inf\{m\geq 0\colon \mathcal{H}^{m}(S)=0\},
\end{align*}
which allows for a more general definition of the dimension of a subset $S\subset \R^{d}$, and agrees with the usual topological definition for smooth manifolds.
\subsubsection{Rectifiable sets}

Rectifiable sets  are a generalisation of smooth manifolds and a central object of study in geometric measure theory, from which one defines currents. 
\begin{definition}
A set $S \subset \R^{d}$ is $m$-rectifiable (for $m\leq d$) if it has Hausdorff dimension $m$, and if there exists bounded $U_{i} \subset \R^{m}$ and Lipschitz functions $f_{i}: U_{i} \longrightarrow \R^{d}$ such that
\begin{align*}
    \mathcal{H}^{m}(S - C) = 0,\quad C\coloneq\cup_{i=1}^{\infty}{f_i(U_{i})}.
\end{align*}
\end{definition}
Unlike charts on smooth manifolds, the cover of the set $S$ is not necessarily smooth, allowing the above definition to encompass piecewise smooth sets such as triangulated surfaces, as well as much rougher objects such as Cantor sets. Such $m$-rectifiable subsets, despite their general definition, can be shown to have well-defined tangent spaces $\mathcal{H}^{m}$-almost everywhere \cite{NCpaper}, which allows one to generalise familiar notions of differential geometry on smooth manifolds in a measure-theoretic sense, such as integration of differential forms over rectifiable sets. This leads to the notion of (rectifiable) currents, as defined in \cref{bg}.

\subsection{Continuity of currents}

The following lemma shows the continuity of the current associated to an oriented rectifiable set.

\begin{lemma}\label[lemma]{simplelemmma}
Let $S \subset \mathbb{R}^d$ be an oriented $m$-rectifiable set with associated current $[S] \in \Omega_0^m(\mathbb{R}^d)^*$ defined by
\[
[S](\omega) \coloneqq \int_S \bigl(\omega(x) \mid \tau_S(x)\bigr)\, d\mathcal{H}^m(x), \qquad \omega \in \Omega_0^m(\mathbb{R}^d),
\]
where $\tau_S(x) = e_1(x) \wedge \cdots \wedge e_m(x)$ is the wedge product of a positively oriented orthonormal basis of the tangent space at $x$. Then,
\[
\bigl|[S](\omega)\bigr| \leq \|\omega\|_\infty\, \mathcal{H}^m(S), \qquad \omega \in \Omega_0^m(\mathbb{R}^d).
\]
\end{lemma}

\begin{proof}
By the standard estimate $\bigl|\int f\bigr| \leq \int |f|$ (triangle inequality for integrals), we have
\[
\bigl|[S](\omega)\bigr| 
= \left|\int_S \bigl(\omega(x) \mid \tau_S(x)\bigr)\, d\mathcal{H}^m(x)\right|
\leq \int_S \bigl|\bigl(\omega(x) \mid \tau_S(x)\bigr)\bigr|\, d\mathcal{H}^m(x).
\]

The pairing $\bigl(\omega(x) \mid \tau_S(x)\bigr)$ is the duality pairing between $\omega(x) \in (\Lambda^m \mathbb{R}^d)^*$ and $\tau_S(x) \in \Lambda^m \mathbb{R}^d$. By definition of the dual norm on $(\Lambda^m \mathbb{R}^d)^*$,
\[
\bigl|\bigl(\omega(x) \mid \tau_S(x)\bigr)\bigr| 
\leq |\omega(x)|_{(\Lambda^m \mathbb{R}^d)^*}\, |\tau_S(x)|_{\Lambda^m \mathbb{R}^d}.
\]
Since $\{e_1(x), \ldots, e_m(x)\}$ is an orthonormal basis of the tangent space, the inner product on $\Lambda^m \mathbb{R}^d$ (cf. \cref{exteriorinnerproddef}) gives
\begin{align*}
|\tau_S(x)|^2 
&= \langle e_1(x) \wedge \cdots \wedge e_m(x),\; e_1(x) \wedge \cdots \wedge e_m(x) \rangle_{\Lambda^m \mathbb{R}^d}\\
&= \det\!\bigl(\langle e_i(x), e_j(x)\rangle\bigr)
= \det(I_m) = 1.
\end{align*}
Therefore, $|\tau_S(x)|_{\Lambda^m \mathbb{R}^d} = 1$, and the pointwise bound simplifies to
\[
\bigl|\bigl(\omega(x) \mid \tau_S(x)\bigr)\bigr| \leq |\omega(x)|_{(\Lambda^m \mathbb{R}^d)^*}.
\]

Substituting the pointwise estimate into the integral,
\[
\bigl|[S](\omega)\bigr| 
\leq \int_S |\omega(x)|_{(\Lambda^m \mathbb{R}^d)^*}\, d\mathcal{H}^m(x)
\leq \sup_{x \in \mathbb{R}^d} |\omega(x)|_{(\Lambda^m \mathbb{R}^d)^*} \cdot \mathcal{H}^m(S)
= \|\omega\|_\infty\, \mathcal{H}^m(S).\qedhere
\]

\end{proof}

\section{Proof of \cref{mainlemma}}\label{APPENDIXSECPROOF}
We now prove \cref{mainlemma}. In the following, we fix an RKHS of forms $W$ as in \cref{newembeddingspace} with spatial kernel $K_{p}$ and $K_{s}=1$. The action of the embedded normal cycle $N(\mathcal{T}) \in W^{*}$ is defined in \cref{newaction}. 

Given a triangulation $\mathcal{T} = \cup_{i=1}^{n}{T_{i}} \subset \R^{3}$ (so $d=3$), a decomposition of $N(\mathcal{T})$ into distinct `planar', `cylindrical' and `spherical' components is derived in \cite{NCpaper} such that
\begin{align}
    N(\mathcal{T}) = N(\mathcal{T})^{\mathrm{pln}} + N(\mathcal{T})^{\mathrm{cyl}} + N(\mathcal{T})^{\mathrm{sph}} \in W^{*},
\end{align}   
and the three components are mutually orthogonal in $W^{*}$. Each component may be further decomposed \cite{NCpaper} into distinct contributions from the triangles and edges of $\mathcal{T}$. We split the proof of \cref{mainlemma} into separate parts for each component, using their decompositions obtained in \cite{NCpaper}.

\subsection{Planar component}
In \cite{NCpaper}, $N(\mathcal{T})^{\mathrm{pln}}$ is derived to be the following sum of currents on the normal bundle:
\begin{align*}
        N(\mathcal{T})^{\mathrm{pln}}  = \sum_{i=1}^{N_{T}}{ [T_{i} \times [ \pm n_{T_{i}}] ] }.
\end{align*}
In the above, $N_{T}$ is the number of triangles of $\mathcal{T}$, $T_{i}$ is the $i$th triangle and $n_{T_{i}}$ a unit normal to the $i$th triangle face. The individual terms in this sum are of the form,
\begin{align*}
   \mathcal{C} = [T \times \{\pm n_{T}\}]
\end{align*}
for some triangle $T$. This component of the normal cycle has the action
\begin{align*}
    \mathcal{C}(\omega) = \int_{T}{\omega_{x}(\tau_{x,n_{T}})\,d\mathcal{H}^{2}(x)} +  \int_{T}{\omega_{x}(\tau_{x,-n_{T}})\,d\mathcal{H}^{2}(x)},\quad \omega \in W.
\end{align*}
Recall that $\tau_{x,u}$ is constructed as a wedge product of an orthonormal basis of the tangent space at $(x,u)$ of the normal bundle. Indeed, one may show \cite{NCpaper} that 
\begin{align*}
    \tau_{x,u} = \begin{pmatrix} e_1(x,u) \\ 0  \end{pmatrix}\wedge  \begin{pmatrix} e_2(x,u) \\ 0  \end{pmatrix}
\end{align*}
where $(e_{1}(x,u),e_{2}(x,u),u)$ is a positively oriented orthonormal basis of $\R^{3}$. Using the identity $ \tau_{x,-u}=- \tau_{x,u}$, the bilinearity and antisymmetry of $\omega_{x}$ as a function of its argument, we have
\begin{align*}
     \mathcal{C}(\omega) &= \int_{T}{(\omega_{x}(\tau_{x,n_{T}}) + \omega_{x}(\tau_{x,-n_{T}}))\,d\mathcal{H}^{2}(x)} \\
     &= \int_{T}{(\omega_{x}(\tau_{x,n_{T}}) + \omega_{x}(-\tau_{x,+n_{T}}))\,d\mathcal{H}^{2}(x)}
   \\  &= \int_{T}{(\omega_{x}(\tau_{x,n_{T}}) - \omega_{x}(\tau_{x,n_{T}})) \,d\mathcal{H}^{2}(x)} = 0,
\end{align*}
for all $\omega \in W$. As a consequence, 
\begin{align*}
    \mathcal{C} = 0 \implies  N(\mathcal{T})^{\mathrm{pln}}=0.
\end{align*}
since $N(\mathcal{T})^{\mathrm{pln}}$ is a sum of such terms. The fact that the planar component vanishes for $K_{s}=1$ is observed in \cite{NCpaper}.
\subsection{Cylindrical component}
\begin{figure}\centering\includegraphics{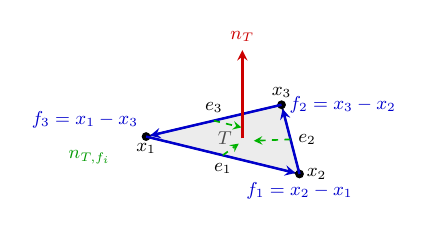}\caption{Notation for edges $e_j$ and vector edges $f_j$ of triangle $T$ with vertices $x_i$. The green dashed lines denote the normals $n_{T,f_i}$.}\label[figure]{fig:triangle}\end{figure}

\begin{remark}\label{notationremark}
  In the following, we consider a triangle $T$ with vertices $x_1,x_2,x_3$ and edges $e_1,e_2,e_3$ connecting vertices $(1,2)$, $(2,3)$, $(3,1)$ respectively. The edges $e_i$ include the endpoints (vertices), with $\tilde e_i$ denoting the open edges without endpoints. The vector edges are defined by $f_1=x_2-x_1$, $f_2=x_3-x_2$, $f_3=x_1-x_3$ (i.e., $f_j$ points away from the vertex $x_j$). The unit normal to $T$ is $n_T\coloneq f_1\times f_2/|f_1\times f_2|$. Note that $f_j\cdot n_T$ is  independent of choice of labels for the vertices and $n_{T,f_{i}}\times f_{i}$ is oriented into the triangle $T$. See \cref{fig:triangle}.
  \end{remark}  

\begin{figure}\centering
  \includegraphics{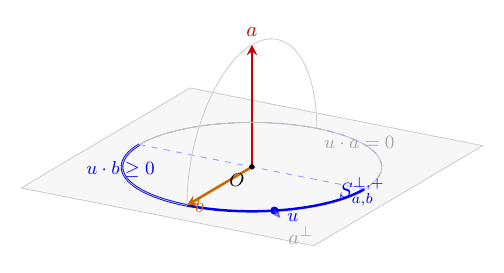}
  \caption{Illustration of $S_{a,b}^{\perp,+}$, the half circle of unit vectors perpendicular to $a$ and with positive inner product with $b$.}\label{halfcirclefig}
  \end{figure}

The cylindrical component of $N(\mathcal{T})$ is given as the following sum of currents on the normal bundle \cite{NCpaper},
\begin{align}\label{NCcylcomp}
      N(\mathcal{T})^{\mathrm{cyl}} 
     = -\sum_{j=1}^{n_{e}}\sum_{\{T\colon e_j\in T\}}{[ e_{j} \times S_{f_{j},-f_{j}\times n_{T}}^{\perp,+} ]}   + \sum_{j=1}^{n_{e}}{[\tilde{e}_{j} \times (e_{j}^{\perp} \cap \mathbb{S}^{2})]},
\end{align}
where $n_{e}$ denotes the total number of edges. The edges $e_{j}$ and open edges $\tilde e_j$ (unoriented) and vector edges $f_j$ are defined in \cref{notationremark}. In the notation of \cite{NCpaper}, one denotes
\begin{equation*}
   v^{\perp} \coloneq \{w \in \R^{3}: v \cdot w = 0\},\quad 
    S_{a,b}^{\perp,+} \coloneq (\mathbb{S}^{2}\cap a^{\perp}) \cap \{u: u\cdot b \geq 0\},
\end{equation*}
so that $S_{a,b}^{\perp,+}$ is an oriented semicircle as in \cref{halfcirclefig}. By definition of $n_T$, $S_{f_j,-f_j\times n_T}^{\perp,+}$ is independent of the choice of labels for the vertices. One may visualise an embedding of $e_{j} \times S_{f_{T},-f_{T}\times n_{T}}^{\perp,+}$ as a half cylinder centered on the edge $e_{j}$ oriented inwards to the triangle. Similiarly, one may visualise $\tilde{e}_{j} \times (e_{j}^{\perp} \cap \mathbb{S}^{2})$ as a full cylinder centered around the edge $e_{j}$. 

  From \cref{NCcylcomp}, the term $N(\mathcal{T})^{\mathrm{cyl}} $ is therefore a sum of cylindrical currents of the form
\begin{equation}\label{cylnc}
    \mathcal{C} = [e_j \times C].
\end{equation}
The sum is over closed edges $e_j$ of triangles $T$ (double counting shared edges) with half circles $C$, and open edges $\tilde e_j$ with full circles $C$. The action of $\mathcal{C}$ on $\omega \in W$ in \cref{cylnc} is given by
\begin{align*}
    \mathcal{C}(\omega) = \int_{e_j \times C}{ \omega_{x}(\tau_{x,u}) \,d\mathcal{H}^{2}(x,u) },\quad \forall\omega \in W.
\end{align*}
As for currents and varifolds and following \cite{NCpaper},  to allow practical metric computation,  we approximate $\mathcal C\approx \tilde{\mathcal C}$ using a midpoint quadrature rule via
\begin{align}\label{cylncapp}
   \hat{\mathcal{C}}(\omega) \coloneq  \lvert e_j \rvert \int_{C}{\omega_{c_j}(\tau_{c_j,u})\,d\mathcal{H}^{1}(u) },
\end{align}
where  $c_j$ is the centre  of edge $e_j$. The remainder $R_j(\omega)\coloneq  \mathcal{C}(\omega) - \hat{\mathcal{C}}(\omega)$   associated to edge $e_j$ is $\mathcal{O}(|e_j|^2)$ and disappears as the length $|e_j|$ of $e_j$  tends to zero.


The $2$-vector $\tau_{x,u}$ is the wedge product of a positively oriented orthonormal basis of the tangent space of $e_j\times C$ at $(x,u)$ and has the following form \cite{NCpaper} on cylindrical components,
\begin{align}\label{taucyl}
    \tau_{x,u} = \begin{pmatrix}  e_1(x,u) \\ 0  \end{pmatrix}\wedge  \begin{pmatrix} 0\\ e_2(x,u) \end{pmatrix}
\end{align}
 where $(e_{1}(x,u),e_{2}(x,u),u)$ is a positively oriented orthonormal basis of $\R^{3}$. At a given point $(x,u)$ on the cylindrical component associated to an edge $f$, one may verify that such a basis may be computed as 
\begin{align*}
(e_{1}(x,u),e_{2}(x,u),u) = \bigg(\frac{f }{\lvert f \rvert},-\frac{f }{\lvert f \rvert} \times u ,u\bigg).
\end{align*}
This expression is independent of $x$, the first component along the cylinder. Substituting this choice of basis into \cref{taucyl,cylncapp} and applying bilinearity of the wedge product yields
\begin{align*}
    \int_{C}{\omega_{c_j}(\tau_{c_j,u}) d\mathcal{H}^{1}(u)}
   & = \int_{C}{\omega_{c_j}\bigg(\begin{pmatrix}  \frac{f_j }{\lvert f_j \rvert}\\ 0  \end{pmatrix}\wedge  \begin{pmatrix} 0\\ -\frac{f_j }{\lvert f_j \rvert} \times u  \end{pmatrix}\bigg) d\mathcal{H}^{1}(u)}\\
&    =\omega_{c_j}\bigg(\begin{pmatrix}  \frac{f_j }{\lvert f_j \rvert} \\ 0  \end{pmatrix}\wedge  \int_{C}{\begin{pmatrix} 0\\ -\frac{f_j }{\lvert f_j \rvert} \times u  \end{pmatrix} d\mathcal{H}^{1}(u)}\bigg)\\
    &    =\omega_{c_j}\bigg(\begin{pmatrix} \frac{f_j }{\lvert f_j \rvert} \\ 0  \end{pmatrix}\wedge  \begin{pmatrix} 0\\ \int_{C}{-\frac{f_j }{\lvert f_j \rvert}\times u 
 \ d\mathcal{H}^{1}(u)} \end{pmatrix} \bigg)\\
 &=\omega_{c_j}\bigg(\begin{pmatrix}  \frac{f_j }{\lvert f_j \rvert} \\ 0  \end{pmatrix}\wedge  \begin{pmatrix} 0\\ -\frac{f_j }{\lvert f_j \rvert} \times \int_{C}{ u 
 \ d\mathcal{H}^{1}(u)} \end{pmatrix}\bigg) 
 = \omega_{c_j}(\gamma_j),
\end{align*}
where one defines
\[ 
    \gamma_j\coloneq \begin{pmatrix}  \frac{f_j }{\lvert f_j \rvert} \\ 0  \end{pmatrix}\wedge  \begin{pmatrix} 0\\ -\frac{f_j }{\lvert f_j \rvert} \times \int_{C}{ u 
 \ d\mathcal{H}^{1}(u)} \end{pmatrix}.
\]

When $C$ is a full circle associated to an edge $e$ (and vector edge $f$), or a half circle associated to an edge $e$ in a triangle $T$, one may compute (see~\cref{halfcircleint})
\begin{align*}
   \int_{C}{ u 
 \ d\mathcal{H}^{1}(u)}  = \begin{cases}
    0, 
    & C=\tilde{e} \times (e^{\perp} \cap \mathbb{S}^{2}),\\
			{2}(-\frac{f }{\lvert f\rvert}\times n_{T}), 
      & C=e \times S_{f,-f\times n_{T}}^{\perp,+},
		 \end{cases}
\end{align*}%
where the first integral evaluates to zero due to symmetry. Upon substitution, this finally yields
\begin{gather*}
    \gamma = \begin{dcases}
      0, & C=\tilde{e}_j \times (e_j^{\perp} \cap \mathbb{S}^{2}),\\ {2}
			\begin{pmatrix}  \frac{f }{\lvert f \rvert} \\ 0  \end{pmatrix}\wedge  \begin{pmatrix} 0\\  -n_{T } \end{pmatrix} , & C=e_j \times S_{f,-f\times n_{T}}^{\perp,+},
		 \end{dcases}
\end{gather*}
where we have used the identity $\frac{f }{\lvert f \rvert} \times (\frac{f}{\lvert f\rvert}\times n_{T})=-n_{T}$ for a vector edge $f$ on the triangle $T$. Substituting into \cref{cylnc,cylncapp}, the above  implies 
\begin{gather*}
 \tilde{\mathcal{C}}=  
\begin{dcases}
    0, & C=\tilde{e}_j \times (e_j^{\perp} \cap \mathbb{S}^{2}),\\ 2
			\begin{pmatrix}  f_j  \\ 0  \end{pmatrix}\wedge  \begin{pmatrix} 0\\  -n_{T } \end{pmatrix} , & C=e_j \times S_{f_j,-f_j\times n_{T}}^{\perp,+},
		 \end{dcases} 
\end{gather*}
so that the full cylinder terms of \cref{NCcylcomp} are zero and only half cylinder terms in the first summation of \cref{NCcylcomp} remain. Note the expression is independent of labelling of the vertices (via definition of  edge vectors $f_j$ and normal $n_T$).

Now we sum over triangles $T$ in \cref{NCcylcomp}:
\begin{align*}
 N(\mathcal{T})^{\mathrm{cyl}}
     = 2\sum_{j=1}^{n_{e}}{\delta_{c_{j}}\sum_{\{T\colon e_j\in \partial T\}}{
			\begin{pmatrix}  f_{j} \\ 0  \end{pmatrix}\wedge  \begin{pmatrix} 0\\  n_{T } \end{pmatrix}} }+\sum_{j=1}^{n_e} R_j,
\end{align*}
where $f_{j}$ is the oriented vector edge of the edge $e_j$. 
We may simplify this further:
\begin{align*}
 N(\mathcal{T})^{\mathrm{cyl}}
     = 2\sum_{j=1}^{n_{e}}{\delta_{c_{j}}\begin{pmatrix}  f_{j} \\ 0  \end{pmatrix} \wedge \bigg(\sum_{\{T\colon e_j\in \partial T\} }{\begin{pmatrix} 0\\  n_{T}  \end{pmatrix}}\bigg) }+\sum_{j=1}^{n_e} R_j\\
     =2\sum_{j=1}^{n_{e}}{\delta_{c_{j}}\begin{pmatrix}  f_{j} \\ 0  \end{pmatrix} \wedge {\begin{pmatrix} 0\\  \sum_{\{T\colon e_j \in \partial T\}} n_{T } \end{pmatrix}} }+\sum_{j=1}^{n_e} R_j.
\end{align*}

This results in a final expression
\begin{align*}
 N(\mathcal{T})^{\mathrm{cyl}}
     = 2\sum_{j=1}^{n_{e}}{\delta_{c_{j}}\gamma_{j} }+\sum_{j=1}^{n_e} R_j,
\end{align*}
where the weights are redefined as
\begin{align}\label{gammadef}
     \gamma_{j}\coloneq 
     \begin{pmatrix}  f_{j} \\ 0  \end{pmatrix} 
     \wedge 
      \begin{pmatrix} 0\\ 
        \sum_{\{T\colon e_j \in \partial T\}}
          n_{T } \end{pmatrix} .
\end{align}
Here, $\gamma_j$ is not required to be a unit 2-vector -- it's a weight in $\Lambda^2(\mathbb{R}^3 \times \mathbb{R}^3)$ arising from integration over the half circles of the normal bundle, and its norm encodes geometric information related to the dihedral angle at the edge. 

\subsection{Spherical component}
Finally, $N(\mathcal{T})^{\mathrm{sph}}$ is shown in \cite{NCpaper} to be equivalent to the normal cycle of the (possibly empty) boundary curve of the triangulation, so $N(\mathcal{T})^{\mathrm{sph}} = N(\partial \mathcal{T})^{\mathrm{sph}}$,
which is further shown to decompose into   
\begin{align}\label{boundarysph}
   N(\mathcal{T})^{\mathrm{sph}} =  N(\partial \mathcal{T})^{\mathrm{sph}} =  \sum_{x_{k} \in \partial \mathcal{T}}{ [\{x_{k}\}\times \mathbb{S}^{2}] } +  \sum_{x_{k} \in \partial \mathcal{T}}\sum_{\{ f_{i}\colon e_k\in \partial f_i\} }{ [\{x_{k}\}\times S_{f_{i}/\lvert f_{i} \rvert}^{+}] }.
\end{align}
where $f_{i}$ are vector edges of the boundary attached to boundary vertices $x_{k}$, and oriented outwards from the vertex. In \cref{boundarysph}, one denotes
\begin{align*}
S_{v}^{+} = \{u \in \mathbb{S}^{2}: u\cdot v \geq 0  \},
\end{align*}
representing oriented half-spheres. Therefore, a generic element of the spherical component of the normal cycle takes the form
\begin{align*}
    \mathcal{C} = [\{x\} \times S],\quad x \in \R^{3}
\end{align*}
for a full or half sphere $S$. Such terms have the following action
\begin{align*}
     \mathcal{C}(\omega) = \int_{S}{\omega_{x}(\tau_{x,u})\,d\mathcal{H}^{2}(u) }\quad \omega \in W 
\end{align*}

On the spherical component $\tau_{x,u}$ has the following generic form \cite{NCpaper} 
\begin{align*}
    \tau_{x,u} = \begin{pmatrix} 0 \\ e_1(x,u)  \end{pmatrix}\wedge  \begin{pmatrix} 0\\ e_2(x,u) \end{pmatrix},
\end{align*}
for appropriate $e_{1}(x,u),e_{2}(x,u)$, which also satisfy $(e_{1}(x,u),e_{2}(x,u),u)$ is a positively oriented orthonormal basis of $\R^{3}$. We may compute an explicit expression for such components by expressing $(e_{1}(x,u),e_{2}(x,u),u)$ in the standard basis of $\R^{3}$. Indeed for fixed $x,u$ one may express (omitting coefficient dependence on $(x,u)$),
\begin{align*}
    e_{1}(x,u) = \sum_{i=1}^{3}{a_{i}e_{i}},\quad e_{2}(x,u) = \sum_{i=1}^{3}{b_{i}e_{i} },\quad u = \sum_{i=1}^{3}{u_{i}e_{i}}.
\end{align*}
By bilinearity and antisymmetry of the wedge product, this substitution in coordinates yields
\begin{align*}
    \tau_{x,u} &= \begin{pmatrix} 0 \\ e_1(x,u)  \end{pmatrix}\wedge  \begin{pmatrix} 0\\ e_2(x,u) \end{pmatrix} = \sum_{i,j=1}^{3}{a_{i}b_{j}\begin{pmatrix} 0 \\ e_i  \end{pmatrix}\wedge  \begin{pmatrix} 0\\ e_j \end{pmatrix}}\\
&    =\sum_{1\leq i<j\leq 3 }{(a_{i}b_{j} - a_{j}b_{i})\begin{pmatrix} 0 \\ e_i  \end{pmatrix}\wedge  \begin{pmatrix} 0\\ e_j \end{pmatrix}  }
\end{align*}
The cyclic coefficients $(a_{i}b_{j} - a_{j}b_{i})$ are nothing but the components of the cross product $e_{1}(x,u) \times e_{2}(x,u)$. Since $(e_{1}(x,u),e_{2}(x,u),u)$ is a positively oriented orthonormal basis, it implies $u = e_{1}(x,u) \times e_{2}(x,u)$. We may thus write
\begin{align*}
  \tau_{x,u}  =  \sum_{i=1}^{3}{u_{i}\tilde{b}_{i}  } 
\end{align*}
in terms of the components of $u$, where $\{\tilde{b}_{i}\}_{i=1}^{3}=\bigg\{\begin{pmatrix} 0 \\ e_i  \end{pmatrix}\wedge  \begin{pmatrix} 0\\ e_j \end{pmatrix}\bigg\}_{{1\leq i<j\leq 3 }}$ are orthonormal. This yields
\begin{align*}
      \omega_{x}(\tau_{x,u}) = \sum_{j=1}^{3}{u_{j}\omega_{x}(\tilde{b}_{j})  } 
\end{align*}
and thus by linearity of $\omega_{x}$ for each $x$,
\begin{align*}
C(\omega)&=\int_{S}{\omega_{x}(\tau_{x,u})\,d\mathcal{H}^{2}(u) } =  \sum_{j=1}^{3}{\bigg(\int_{S}u_{j}\,d\mathcal{H}^{2}(u))\bigg)\omega_{x}(\tilde{b}_{j}) } \\&=  {\omega_{x}\bigg(\sum_{j=1}^{3}\bigg(\int_{S}u_{j}\,d\mathcal{H}^{2}(u)\bigg)\tilde{b}_{j}\bigg) }.
\end{align*}
Thus, one may represent elements of the spherical component as
\begin{align*}
    \mathcal{C} = \delta_{x}\alpha,\qquad \alpha\coloneq \sum_{j=1}^{3}\bigg(\int_{S}u_{j}\,d\mathcal{H}^{2}(u)\bigg)\tilde{b}_{j}.
\end{align*}
Evaluating the spherical unit normal integral  (see \cref{hemisphereint}) when $S$ is a half-sphere or a full sphere, we see
\begin{gather*}
    \int_{S}{u\,d\mathcal{H}^{2}(u) }  = \begin{dcases}
        \pi v, & S=S_{v}^{+},\\ 0, & S=\mathbb{S}^{2},
		 \end{dcases} 
\end{gather*}
so that the full spherical terms vanish due to symmetry, and
\begin{gather*}
    \mathcal{C} = [\{x\}\times S] = \begin{dcases}
        \delta_{x}\alpha ,\quad \alpha\coloneq \pi \sum_{j=1}^{3}v_{j}\tilde{b}_{j}, & S=S_{v}^{+},\\ 0 , & S=\mathbb{S}^{2}.
		 \end{dcases} 
\end{gather*}
 Substituting the above expressions into the half-sphere terms of \cref{boundarysph} (where $v$ are normalised edges) yields
\[
N(\mathcal{T})^{\mathrm{sph}} = \pi\sum_{x_{k} \in \partial \mathcal{T }}\delta_{x_{k}} \sum_{\{f_{i}\colon x_k\in \partial f_i\}}\sum_{j=1}^{3}{\frac{f_{ij}}{\lvert f_{i} \rvert} }\tilde{b}_{j}, 
\]
where $f_{ij}$ is the $j$th component of $f_i$. Then,
\begin{align*}
  N(\mathcal{T})^{\mathrm{sph}}  =  \pi\sum_{x_{k} \in \partial \mathcal{T }}\delta_{x_{k}}{\alpha_{k}} ,
\end{align*}
where one defines vertex weights
\begin{align}\label{alphweights}
\alpha_{k} \coloneq \sum_{\{f_{i}\colon x_k\in \partial f_i\}}\sum_{j=1}^{3}{\frac{f_{ij}}{\lvert f_{i} \rvert} }\tilde{b}_{j} \in \Lambda^{2}(\R^{3}\times \R^{3}),
\end{align}
where we recall $f_{i}$ are oriented outward from $x_{k}$.
\subsection{Full decomposition and inner product}
Combining the expressions derived so far gives the following Dirac delta decomposition of the normal cycle:
\begin{align*}
    N(\mathcal{T}) = N(\mathcal{T})^{\mathrm{cyl}} + N(\mathcal{T})^{\mathrm{sph}} = 2 \sum_{i=1}^{n_{e}}{\delta_{c_{i}}\gamma_{i}} +  \pi\sum_{x_{k} \in \partial \mathcal{T }}\delta_{x_{k}}{\alpha_{k}}
\end{align*}
for weights defined in \cref{alphweights,gammadef}. The first term is the cylindrical component $N(\mathcal{T})^{\mathrm{cyl}}$, and the second term is the spherical component $N(\mathcal{T})^{\mathrm{sph}}$.

One may verify that this Dirac delta decomposition induces the inner product \cref{NC_disc_met} between two normal cycles $N(\mathcal{T}),N(\mathcal{T}')$. Indeed, denoting $\alpha_{i},\gamma_{i}$ weights of $ N(\mathcal{T})$, and $\alpha_{i}',\gamma_{i}'$ weights of $N(\mathcal{T}')$, one observes
\begin{gather}\begin{split}
   & \langle N(\mathcal{T}),N(\mathcal{T}') \rangle_{W^{*}} = \langle N(\mathcal{T})^{\mathrm{cyl}},N(\mathcal{T}')^{\mathrm{cyl}}\rangle_{W^{*}} + \langle N(\mathcal{T})^{\mathrm{sph}} ,N(\mathcal{T}')^{\mathrm{sph}}\rangle_{W^{*}} \\
   &= 4\sum_{i=1}^{n_{e}}\sum_{j=1}^{n_{e}'}{\langle \delta_{c_{i}}\gamma_{i},\delta_{c_{j}'}\gamma_{j}' \rangle_{W^{*}}} + \pi^2\sum_{x_{k} \in \partial\mathcal{T }}\sum_{y_{l} \in \partial\mathcal{T }'}{\langle \delta_{x_{k}}\alpha_{k},\delta_{y_{l}}\alpha_{l}' \rangle_{W^{*}} }\\
    &= 4 \sum_{i=1}^{n_{e}}\sum_{j=1}^{n_{e}'}{K_{p}(c_{i},c_{j}')\langle\gamma_{i},\gamma_{j}' \rangle_{\Lambda^{2}(\R^{3}\times \R^{3})}} + \pi^2\sum_{x_{k} \in \partial\mathcal{T }}\sum_{y_{l} \in \partial\mathcal{T }'}{K_p(x_{k},y_{l})\langle \alpha_{k},\alpha_{l}' \rangle_{\Lambda^{2}(\R^{3}\times \R^{3})} },
\end{split}\label{FINALINE}\end{gather}
where we use the orthogonality of spherical and cylindrical terms \cite{NCpaper} in the first line, and the definition of RKHS inner product in the final line. Using the definition in \cref{exteriorinnerproddef} of the inner product on $\Lambda^{2}(\R^{3}\times \R^{3})$, one observes for $a,b,c,d \in \R^{3}$,
\begin{align*}
      \bigg\langle \begin{pmatrix} a \\ 0  \end{pmatrix}\wedge \begin{pmatrix} 0 \\ b  \end{pmatrix} ,\begin{pmatrix} c \\ 0  \end{pmatrix} \wedge \begin{pmatrix} 0 \\ d  \end{pmatrix}   \bigg\rangle_{\Lambda^{2}(\R^{3}\times \R^{3})}& = (a\cdot c)(b\cdot d),\\
    \bigg\langle \begin{pmatrix} 0 \\ a  \end{pmatrix}\wedge \begin{pmatrix} 0 \\ b  \end{pmatrix},\begin{pmatrix} 0 \\ c  \end{pmatrix} \wedge \begin{pmatrix} 0 \\ d  \end{pmatrix}  \bigg\rangle_{\Lambda^{2}(\R^{3}\times \R^{3})}&=(a\cdot c)(b\cdot d) - (b\cdot c)(a\cdot d). 
\end{align*}
Substituting for the inner product of weights in \cref{FINALINE}, this yields the following for the spherical weights (where we denote $x_{k},y_{l}$ boundary vertices of $\mathcal{T},\mathcal{T}'$, respectively),
\begin{align*}
    \langle \alpha_{k},\alpha_{l}' \rangle_{\Lambda^{2}(\R^{3}\times \R^{3})}
    &= \bigg\langle \sum_{\{f_{i}\colon x_k\in \partial f_i\}}\sum_{p=1}^{3}{\frac{f_{ip}}{\lvert f_{i} \rvert} }\tilde{b}_{p}, \sum_{\{g_{j}\colon y_l\in \partial g_j\}}\sum_{q=1}^{3}{\frac{g_{jq}}{\lvert g_{j} \rvert} }\tilde{b}_{q}\bigg\rangle_{\Lambda^{2}(\R^{3}\times \R^{3})} \\
    &= \sum_{\{f_i\colon x_k\in \partial f_i\}}
    \sum_{\{g_j\colon  y_{l}\in \partial g_j\}}\sum_{r=1}^{3}{\frac{f_{ir}}{\lvert f_{i} \rvert} }{\frac{g_{jr}}{\lvert g_{j} \rvert} }\langle \tilde{b}_{r},\tilde{b}_{r}\rangle_{\Lambda^{2}(\R^{3}\times \R^{3})}\\
&=    \sum_{\{f_i\colon x_k\in \partial f_i\}}
    \sum_{\{g_j\colon  y_{l}\in \partial g_j\}}{\bigg\langle \frac{f_{i}}{\lvert f_{i} \rvert},\frac{g_{j}}{\lvert g_{j} \rvert} \bigg\rangle_{\R^{3}}}\\&=   \bigg\langle \sum_{\{f_{i}\colon x_k \in \partial f_i\}}{ \frac{f_{i}}{\lvert f_{i} \rvert} },\sum_{\{g_{j}\colon y_l\in \partial g_j\}}{\frac{g_{j}}{\lvert g_{j} \rvert}} \bigg\rangle_{\R^{3}},
\end{align*}
where we use the orthonormality of $\tilde{b}_{j} \in \Lambda^{2}(\R^{3}\times \R^{3})$ in the second line. Similarly, the following holds for cylindrical weights (where we denote $f_{i},g_{j}$ edges in $\mathcal{T},\mathcal{T}'$, respectively),
\begin{align*}
    &\langle\gamma_{i},\gamma_{j}' \rangle_{\Lambda^{2}(\R^{3}\times \R^{3})}  \\
    &=
    \bigg\langle\begin{pmatrix}  f_{i}\\ 0  \end{pmatrix}
    \wedge { \begin{pmatrix} 0\\  \underset{{\{T\colon f_i \in \partial T\} }}{\sum}n_{T,f_{i} } \end{pmatrix} }, \begin{pmatrix}  g_{j} \\ 0  \end{pmatrix} \wedge { \begin{pmatrix} 0\\   \underset{\{T'\colon g_j\in \partial T'\} }{\sum}n_{T',g_{j} } \end{pmatrix} } \bigg\rangle_{\Lambda^{2}(\R^{3}\times \R^{3})}  \\
 &   = \langle f_{i},g_{j} \rangle_{\R^{3}}\bigg\langle \sum_{\{T\colon f_i \in \partial T\}} n_{T,f_{i} },\sum_{\{T'\colon g_j\in \partial T'\} }n_{T',g_{j} }  \bigg\rangle_{\R^{3}}.
\end{align*}
Substituting the identities above into \cref{FINALINE} yields \cref{NC_disc_met}:
\begin{align*}
     &\langle N(\mathcal{T}),N(\mathcal{T}') \rangle_{W^{*}}  \\&= 4\sum_{i=1}^{n_{e}}\sum_{j=1}^{n_{e}'}{K_{p}(c_{i},c_{j}')\langle f_{i},g_{j} \rangle_{\R^{3}}\bigg\langle \sum_{\{T\colon f_i \in\partial T\}}n_{T,f_{i} },\sum_{\{T'\colon g_j\in\partial T'\}}n_{T',g_{j} } \bigg\rangle_{\R^{3}}} \\&\quad + \
     \pi^{2}\sum_{x_{k} \in \partial\mathcal{T }}\sum_{y_{l} \in \partial\mathcal{T }'}{K_p(x_{k},y_{l}) \bigg\langle \sum_{\{f_{i} \colon x_k \in \partial f_i\}}{ \frac{f_{i} }{\lvert f_{i} \rvert} },\sum_{\{g_{j}\colon y_l\in \partial g_j\}}{\frac{g_{j} }{\lvert g_{j} \rvert}} \bigg\rangle_{\R^{3}} }.
\end{align*}
\subsection{Half-circle and hemi-sphere integral}
We give two integral calculations for completeness.
\begin{lemma}\label[lemma]{halfcircleint}
  Let $e$ be an edge of a triangle $T$ with vector edge $f = x_2 - x_1$ and unit normal $n_T$.
  For the half-circle,
  \[
    C = e \times S^{\perp,+}_{f,\,-f\times n_T}
    \;=\;
    \bigl\{\, u \in \mathbb{S}^2 \colon
      u \cdot f  = 0,\;
       u\cdot (-f \times n_T) \geq 0
    \,\bigr\},
  \]
  one has
  \[
    \int_C u \, d\mathcal{H}^1(u)
    \;=\;
    2\!\left(-\frac{f}{|f|} \times n_T\right).
  \]
  \end{lemma}
  
  \begin{proof} Denote the unit vector $\hat{f} = {f}/{|f|}$.
  Since $n_T \perp \hat{f}$, the pair $\{n_T,\, \hat{f} \times n_T\}$ is an
  orthonormal basis of $f^\perp \cap \mathbb{R}^3$.
  Parameterise the great circle $\mathbb{S}^2 \cap f^\perp$ by arc length $\theta$ so 
  \[
    u(\theta)
    = \cos\theta\; n_T + \sin\theta\;(\hat{f} \times n_T),
    \qquad \theta \in [-\pi, \pi),
    \qquad d\mathcal{H}^1 = d\theta.
  \]
  The half-circle condition $ u(\theta)\cdot( -f \times n_T) \geq 0$
  reads $-|f|\sin\theta \geq 0$ (using $n_T\cdot (f\times n_T)=0$ and $(\hat f\times n_T)\cdot(f\times n_T)=|f|$), which is equivalent to $\theta \in [-\pi, 0]$.
  Therefore,
  \[
    \begin{aligned}
    \int_C u \, d\mathcal{H}^1(u)
    &= \int_{-\pi}^{0}
         \Bigl(\cos\theta\; n_T + \sin\theta\;(\hat{f} \times n_T)\Bigr)
       \, d\theta \\
    &= n_T \int_{-\pi}^{0} \cos\theta \, d\theta
      \;+\;
      (\hat{f} \times n_T) \int_{-\pi}^{0} \sin\theta \, d\theta\\
 &   = 0 \cdot n_T + (-2)(\hat{f} \times n_T)\\
  &  = 2\!\left(-\frac{f}{|f|} \times n_T\right). \qedhere
  \end{aligned}
  \]
  \end{proof}

  \begin{lemma}\label[lemma]{hemisphereint}
    Let $v \in \mathbb{S}^2$ be a unit vector and
    $S^+_v = \{ u \in \mathbb{S}^2 \colon  u\cdot v  \geq 0 \}$
    the closed hemisphere oriented by $v$. Then
    \[
      \int_{S^+_v} u \, d\mathcal{H}^2(u) = \pi v.
    \]
    \end{lemma}
    
    \begin{proof}
    By rotational symmetry it suffices to take $v = e_3 = (0,0,1)$.
    In spherical coordinates,
    \[
      u(\theta,\phi) = (\sin\theta\cos\phi,\; \sin\theta\sin\phi,\; \cos\theta),
      \quad \theta\in[0,\tfrac{\pi}{2}],\ \phi\in[0,2\pi),
    \]
    with $d\mathcal{H}^2 = \sin\theta\, d\theta\, d\phi$. The $x$ and $y$
    components of $\int_{S^+_{e_3}} u\, d\mathcal{H}^2$ vanish on integrating
    out $\phi$, while the $z$ component is
    $\int_0^{2\pi}\!\int_0^{\pi/2} \cos\theta\sin\theta\, d\theta\, d\phi = \pi$.
    Hence $\int_{S^+_{e_3}} u\, d\mathcal{H}^2 = \pi e_3$, and the general case follows by rotation.
    \end{proof}

\section*{Declarations}

\bmhead{Funding}
Allen Paul was funded by the EPSRC Centre for Doctoral Training in Statistical Applied Mathematics (SAMBa) at the University of Bath. The other authors received no specific funding for this work.

\bmhead{Consent to Publish}
Not applicable.

\bmhead{Ethics and Consent to Participate}
Not applicable.

\bibliography{references}

@article{NCpaper,
author = {Roussillon, Pierre and Glaunès, Joan Alexis},
address = {New York},
copyright = {Springer Science+Business Media, LLC, part of Springer Nature 2019},
issn = {0924-9907},
journal = {Journal of Mathematical Imaging and Vision},
language = {eng},
number = {8},
pages = {1069-1095},
publisher = {Springer US},
title = {Representation of Surfaces with Normal Cycles and Application to Surface Registration},
volume = {61},
year = {2019},
doi={10.1007/s10851-019-00888-x}
}

@article{Federer1959,
  author    = {Federer, Herbert},
  title     = {Curvature Measures},
  journal   = {Transactions of the American Mathematical Society},
  volume    = {93},
  number    = {3},
  pages     = {418--491},
  year      = {1959},
  doi       = {10.1090/S0002-9947-1959-0110078-1}
}

@phdthesis{NCThesis,
  TITLE = {{Normal cycle models for deformation analysis}},
  AUTHOR = {Roussillon, Pierre},
  URL = {https://theses.hal.science/tel-02180601},
  NUMBER = {2017USPCB073},
  SCHOOL = {{Universit{\'e} Sorbonne Paris Cit{\'e}}},
  YEAR = {2017},
  MONTH = Nov,
  TYPE = {Theses},
  HAL_ID = {tel-02180601},
  HAL_VERSION = {v1},
}

@article{zahle1986integral,
  title={Integral and current representation of {F}ederer's curvature measures},
  author={Z{\"a}hle, Martina},
  journal={Archiv der Mathematik},
  volume={46},
  number={6},
  pages={557--567},
  year={1986},
  publisher={Springer},
  doi={10.1007/bf01195026}
}

@article{fu1989monge,
  title={Monge--{A}mp{\`e}re functions, I},
  author={Fu, Joseph HG},
  journal={Indiana University Mathematics Journal},
  volume={38},
  number={3},
  pages={745--771},
  year={1989},
  publisher={JSTOR},
  doi={10.1512/iumj.1989.38.38035}
}

@inproceedings{CohenSteiner2003,
  author    = {Cohen-Steiner, David and Morvan, Jean-Marie},
  title     = {Restricted {D}elaunay Triangulations and Normal Cycle},
  booktitle = {Proceedings of the Nineteenth Annual Symposium on Computational Geometry (SCG '03)},
  pages     = {312--321},
  year      = {2003},
  publisher = {ACM},
  address   = {New York, NY, USA},
  doi       = {10.1145/777792.777839}
}

@incollection{wintgen1982normal,
  author    = {Wintgen, P.},
  title     = {Normal cycle and integral curvature for polyhedra in {R}iemannian manifolds},
  booktitle = {Differential Geometry},
  editor    = {Soos, Gy. and Szenthe, J.},
  publisher = {North-Holland},
  address   = {Amsterdam},
  year      = {1982},
}

@book{Morvan2008,
  author    = {Morvan, Jean-Marie},
  title     = {Generalized Curvatures},
  series    = {Geometry and Computing},
  volume    = {2},
  publisher = {Springer},
  address   = {Berlin, Heidelberg},
  year      = {2008},
  doi       = {10.1007/978-3-540-73792-6}
}

@article{SparseNystromCurrVar,
      title={Compression  of Currents and Varifolds}, 
  journal = {SIAM Journal on Imaging Sciences},
  volume = {19},
  number = {1},
  pages = {327-363},
  year = {2026},
  doi = {10.1137/24M1699656}
}

@inproceedings{Kaltenmark2017AGF,
  title={A General Framework for Curve and Surface Comparison and Registration with Oriented Varifolds},
  author={Irene Kaltenmark and Benjamin Charlier and Nicolas Charon},
  journal={2017 IEEE Conference on Computer Vision and Pattern Recognition (CVPR)},
  year={2017},
  pages={4580-4589},
  url={https://openaccess.thecvf.com/content_cvpr_2017/papers/Kaltenmark_A_General_Framework_CVPR_2017_paper.pdf}
}

@article{Jakob2015IFA,
  author    = {Jakob, Wenzel and Tarini, Marco and Panozzo, Daniele and Sorkine-Hornung, Olga},
  title     = {Instant Field-Aligned Meshes},
  journal   = {ACM Transactions on Graphics (Proceedings of SIGGRAPH Asia)},
  volume    = {34},
  number    = {6},
  articleno = {189},
  pages     = {1--15},
  year      = {2015},
  month     = nov,
  doi       = {10.1145/2816795.2818078},
  publisher = {ACM}
}

@article{Garland1997QEM,
  author    = {Garland, Michael and Heckbert, Paul S.},
  title     = {Surface Simplification using Quadric Error Metrics},
  journal   = {Proceedings of SIGGRAPH 97},
  pages     = {209--216},
  year      = {1997},
  doi       = {10.1145/258734.258849}
}

@article{Pymeshlab,
  author    = {Muntoni, Alessandro and Cignoni, Paolo},
  title     = {{PyMeshLab}},
  year      = {2021},
  doi       = {10.5281/zenodo.4438750},
  publisher = {Zenodo}
}

@article{RKHS,
doi={10.2307/1990404},
  author = {Aronszajn, N.},
  description = {CiteULike: Theory of reproducing kernels},
  journal = {Transactions of the American Mathematical Society},
  number = 3,
  pages = {337--404},
  title = {Theory of {R}eproducing {K}ernels},
  url = {http://dx.doi.org/10.2307/1990404},
  volume = 68,
  year = 1950
}

@article{DURRLEMAN,
title = {Statistical models of sets of curves and surfaces based on currents},
journal = {Medical Image Analysis},
volume = {13},
number = {5},
pages = {793-808},
year = {2009},
issn = {1361-8415},
doi = {https://doi.org/10.1016/j.media.2009.07.007},
url = {https://www.sciencedirect.com/science/article/pii/S1361841509000620},
author = {Stanley Durrleman and Xavier Pennec and Alain Trouvé and Nicholas Ayache}
}

@INPROCEEDINGS{DAC,
  author={Cherfaoui, Farah and Kadri, Hachem and Ralaivola, Liva},
  booktitle={ICASSP 2022 - 2022 IEEE International Conference on Acoustics, Speech and Signal Processing (ICASSP)}, 
  title={Scalable Ridge Leverage Score Sampling for the {N}yström Method}, 
  year={2022},
  volume={},
  number={},
  pages={4163-4167},
  doi={10.1109/ICASSP43922.2022.9747039}}

@article{Charon,
issn = {1936-4954},
journal = {SIAM Journal on Imaging Sciences},
pages = {2547--2580},
volume = {6},
publisher = {Siam Publications},
number = {4},
year = {2013},
title = {The varifold representation of nonoriented shapes for diffeomorphic registration},
copyright = {Copyright 2014 Elsevier B.V., All rights reserved.},
language = {eng},
address = {PHILADELPHIA},
author = {Charon, Nicolas and Trouvé, Alain},
doi={10.1137/130918885},
}

@InProceedings{GlaunesCurrents,
author="Vaillant, Marc
and Glaun{\`e}s, Joan",
editor="Christensen, Gary E.
and Sonka, Milan",
title="Surface Matching via Currents",
booktitle="Information Processing in Medical Imaging",
year="2005",
publisher="Springer Berlin Heidelberg",
address="Berlin, Heidelberg",
pages="381--392",
doi={10.1007/11505730_32},
isbn="978-3-540-31676-3"
}

@article{VarCompression,
issn = {1615-3375},
journal = {Foundations of Computational Mathematics},
pages = {1317--1361},
volume = {21},
publisher = {Springer US},
number = {5},
doi={DOI: 10.1007/s10208-020-09484-7},
year = {2021},
title = {Metrics, Quantization and Registration in Varifold Spaces},
copyright = {SFoCM 2021},
language = {eng},
address = {New York},
author = {Hsieh, Hsi-Wei and Charon, Nicolas},
}

@book{Younes,
series = {Applied Mathematical Sciences},
issn = {0066-5452},
volume = {171},
publisher = {Springer Berlin / Heidelberg},
isbn = {9783642263484},
year = {2010},
title = {Shapes and Diffeomorphisms},
copyright = {Springer-Verlag Berlin Heidelberg 2010},
language = {eng},
address = {Berlin, Heidelberg},
author = {Younes, Laurent},
doi={10.1007/978-3-642-12055-8} 
}

@article{Keops,
  author  = {Benjamin Charlier and Jean Feydy and Joan Alexis Glaunès and François-David Collin and Ghislain Durif},
  title   = {Kernel Operations on the {GPU}, with Autodiff, without Memory Overflows},
  journal = {Journal of Machine Learning Research},
  year    = {2021},
  volume  = {22},
  number  = {74},
  pages   = {1-6},
doi={10.48550/arXiv.2004.11127}

}

@inproceedings{Musco,
 author = {Musco, Cameron and Musco, Christopher},
 booktitle = {Advances in Neural Information Processing Systems},
 editor = {I. Guyon and U. Von Luxburg and S. Bengio and H. Wallach and R. Fergus and S. Vishwanathan and R. Garnett},
 pages = {},
 publisher = {Curran Associates, Inc.},
 title = {Recursive Sampling for the {N}ystrom Method},
 url = {https://proceedings.neurips.cc/paper_files/paper/2017/file/a03fa30821986dff10fc66647c84c9c3-Paper.pdf},
 volume = {30},
address={Long Beach},
 year = {2017}
}

@incollection{FidelityMetrics,
title = {12 - Fidelity metrics between curves and surfaces: currents, varifolds, and normal cycles},
editor = {Xavier Pennec and Stefan Sommer and Tom Fletcher},
booktitle = {Riemannian Geometric Statistics in Medical Image Analysis},
publisher = {Academic Press},
pages = {441-477},
address={London},
year = {2020},
isbn = {978-0-12-814725-2},
doi = {https://doi.org/10.1016/B978-0-12-814725-2.00021-2},

author = {Nicolas Charon and Benjamin Charlier and Joan Glaunès and Pietro Gori and Pierre Roussillon}
}

@article{Sparseadaptiveparam,
issn = {0920-5691},
journal = {International journal of computer vision},
pages = {161--183},
volume = {101},
publisher = {Springer US},
number = {1},
year = {2013},
title = {Sparse Adaptive Parameterization of Variability in Image Ensembles},
copyright = {Springer Science+Business Media, LLC 2012},
language = {eng},
address = {Boston},
author = {Durrleman, Stanley and Allassonnière, Stéphanie and Joshi, Sarang},
doi={10.1007/s11263-012-0556-1},
}

@book{Federer,
author = {Federer, Herbert},
address = {Berlin, Heidelberg},
booktitle = {Geometric measure theory},
doi={10.1007/978-3-642-62010-2},
isbn = {3-540-04505-8},
language = {eng},
publisher = {Springer-Verlag},
series = {Classics in Mathematics},
volume=153,
title = {Geometric measure theory },
year = {1996},
}

@book{FrankMorgan,
author = {Morgan, Frank},
address = {Amsterdam},
booktitle = {Geometric measure theory: a beginner's guide},
edition = {Fifth},
doi={10.1016/C2015-0-01918-9},
isbn = {9780128045275},
language = {eng},
publisher = {Academic Press},
title = {Geometric measure theory: a beginner's guide },
year = {2016},
}

@inproceedings{Chamferref,
author = {Wu, Tong and Pan, Liang and Zhang, Junzhe and Wang, Tai and Liu, Ziwei and Lin, Dahua},
title = {Density-aware {C}hamfer distance as a comprehensive metric for point cloud completion},
year = {2021},
isbn = {9781713845393},
publisher = {Curran Associates Inc.},
address = {Red Hook, NY, USA},
booktitle = {Proceedings of the 35th International Conference on Neural Information Processing Systems},
articleno = {2228},
numpages = {13},
series = {NIPS 2021},
url={https://doi.org/10.48550/arXiv.2111.12702}
}

@inproceedings{Hausdorffref,
author = {Dubuisson, M.-P. and Jain, A.K.},
booktitle = {Proceedings of 12th International Conference on Pattern Recognition},
isbn = {0818662654},
language = {eng},
pages = {566-568},
doi={https://doi.org/10.1109/ICPR.1994.576361},
volume=1,
publisher = {IEEE},
address="Israel",

title = {A modified {H}ausdorff distance for object matching},
volume = {1},
year = {1994},
}

@article{FINDINGCORRESPONDENCE,
author = {van Kaick, Oliver and Zhang, Hao and Hamarneh, Ghassan and Cohen-Or, Daniel},
address = {Oxford, UK},
copyright = {2011 The Authors Computer Graphics Forum © 2011 The Eurographics Association and Blackwell Publishing Ltd.},
issn = {0167-7055},
journal = {Computer graphics forum},
language = {eng},
doi={10.1111/j.1467-8659.2011.01884.x},
number = {6},
pages = {1681-1707},
publisher = {Blackwell Publishing Ltd},
title = {A Survey on Shape Correspondence},
volume = {30},
year = {2011},
}

@book{Sternberg,
author = {Sternberg, Shlomo.},
address = {New York, N.Y},
booktitle = {Lectures on differential geometry},
edition = {2nd},
isbn = {0828403163},
language = {eng},
publisher = {Chelsea},
title = {Lectures on differential geometry },
year = {1983},
}

@article{ANIMdataset,
author = {Sumner, Robert W. and Popović, Jovan},
title = {Deformation transfer for triangle meshes},
volume = {23},
year = {2004},
address = {New York, NY, USA},
copyright = {ACM},
issn = {0730-0301},
journal = {ACM transactions on graphics},
language = {eng},
number = {3},
pages = {399-405},
publisher = {ACM},
doi={10.1145/1015706.1015736}
}

@article{Thingi10K,
  title={Thingi10K: A Dataset of 10,000 3D-Printing Models},
  author={Zhou, Qingnan and Jacobson, Alec},
  journal={arXiv preprint arXiv:1605.04797},
  year={2016},
doi={10.48550/arXiv.1605.04797}
}

@inproceedings{ARCFlow,
  author    = {Hartshorne, Adam and Paul, Allen and Shardlow, Tony and Campbell, Neill D.~F.},
  title     = {{ARC-Flow}: Articulated, Resolution-Agnostic, Correspondence-Free Matching and Interpolation of {3D} Shapes Under Flow Fields},
  year      = {2025},
  eprint    = {2503.02606},
  archivePrefix = {arXiv},
  primaryClass  = {cs.CV},
  url       = {https://arxiv.org/abs/2503.02606}
}

@article{CurvatureEnthusiasm,
  author    = {Hartshorne, Adam and Paul, Allen and Shardlow, Tony and Campbell, Neill D.~F.},
  title     = {Curvature Enthusiasm: Correspondence-Free Interpolation and Matching of Articulated {3D} Shapes using Compressed Normal Cycles},
  journal   = {ACM Transactions on Graphics},
  volume    = {44},
  number    = {6},
  articleno = {230},
  year      = {2025},
  month     = dec,
  pages     = {1--25},
  doi       = {10.1145/3763366}
}
\end{document}